\documentclass[12pt]{article}
\usepackage{a4wide}
\usepackage[centertags,leqno]{amsmath}
\usepackage{amssymb,amsfonts}
\usepackage{color}
\usepackage{graphicx}
\usepackage[colorlinks=true]{hyperref}

\newtheorem{theorem}{Theorem}[section]

\newtheorem{corollary}[theorem]{Corollary}
\newtheorem{example}[theorem]{Example}

\newtheorem{lemma}[theorem]{Lemma}
\newtheorem{proposition}[theorem]{Proposition}
\newtheorem{remark}[theorem]{Remark}

\definecolor{violet}{rgb}{0.5,0,0.5}
\definecolor{orange}{cmyk}{0,0.3,0.7,0}

\numberwithin{equation}{section}

\newcommand{\qed}{\rule{2mm}{2mm}}

\newcommand{\eqdef}{\stackrel{{\mathrm {def}}}{=}}
\newcommand{\eps}{\varepsilon}
\renewcommand{\colon}{:\,}
\newcommand{\proof}{{\em Proof. }}

\newcommand{\RR}{\mathbb{R}}

\newcommand{\NN}{\mathbb{N}}
\newcommand{\ZZ}{\mathbb{Z}}

\newcommand{\ee}{\mathrm{e}}

\newcommand{\wstarconverge}{\stackrel{*}{\rightharpoonup}}
\newcommand{\esssup}{\mathop{\rm {ess\,sup}}\limits} 
\newcommand{\essinf}{\mathop{\rm {ess\,inf}}\limits}

\begin{document}

\begingroup\Large\bf
\begin{center}
  Convergence to travelling waves\\
  in the Fisher\--Kolmogorov equation\\
  with a non\--Lipschitzian reaction term
\end{center}
\endgroup

\vspace{0.5cm}
\begingroup\rm
\begin{center}
\centerline{Pavel {\sc Dr\'{a}bek}
}
\begin{tabular}{c}
  Department of Mathematics and\\
  N.T.I.S. (Center of New Technologies for Information Society)\\
  University of West Bohemia\\
  P.O.~Box 314\\
  CZ-306 14 Plze\v{n}, Czech Republic\\
  {\it e}-mail: {\tt pdrabek@kma.zcu.cz}
\end{tabular}
\end{center}
\endgroup

\vspace{-0.3cm}
\begin{center}
and
\end{center}

\vspace{-0.3cm}
\begingroup\rm
\begin{center}
\centerline{Peter {\sc Tak\'{a}\v{c}}
}
\begin{tabular}{c}
  Institut f\"ur Mathematik\\
  Universit\"at Rostock\\
  Ulmenstra{\ss}e~69, Haus~3\\
  D-18055 Rostock, Germany\\
  {\it e}-mail: {\tt peter.takac@uni-rostock.de}
\end{tabular}
\end{center}
\endgroup

\vspace{-0.1cm}

\begin{center}
\today
\end{center}

\vspace{0.1cm}


\baselineskip=12pt
\noindent
\begingroup\footnotesize
{\bf {\sc Abstract.}}
We consider the semilinear
Fisher\--Kolmogorov\--Petrovski\--Piscounov equation for
the advance of an advantageous gene in biology.
Its nonsmooth reaction function $f(u)$ allows for
the introduction of travelling waves with a new profile.
We study existence, uniqueness, and long\--time asymptotic behavior
of the solutions $u(x,t)$, $(x,t)\in \RR\times \RR_+$.
We prove also the existence and uniqueness (up to a spatial shift)
of a travelling wave $U$.
Our main result is the uniform convergence (for $x\in \RR$)
of every solution $u(x,t)$ of the Cauchy problem to
a single travelling wave $U(x-ct + \zeta)$ as $t\to \infty$.
The speed $c$ and the travelling wave $U$ are determined uniquely by $f$,
whereas the shift $\zeta$ is determined by the initial data.
\endgroup

\vfill
\par\vspace*{0.2cm}
\noindent
\begin{tabular}{ll}
{\bf Running head:}
& Convergence to a travelling wave in the FKPP equation\\
\end{tabular}

\par\vspace*{0.2cm}
\noindent
\begin{tabular}{ll}
{\bf Keywords:}
& Fisher\--Kolmogorov equation, nonsmooth reaction function,\\
& travelling waves, solutions of the Cauchy problem,\\
& long\--time behavior\\
\end{tabular}

\par\vspace*{0.2cm}
\noindent
\begin{tabular}{ll}
{\bf 2010 Mathematics Subject Classification:}
& Primary   35Q92, 35K91;\\
& Secondary 35K58, 92B05.\\
\end{tabular}
 
\newpage

\baselineskip=14pt
\section{Introduction}
\label{s:Intro}

We are concerned with the long\--time asymptotic behavior of solutions
to the Cauchy problem for the favorite {\it Fisher\--KPP equation\/}
(or {\it Fisher\--Kolmogorov equation\/}):
\begin{equation}
\label{e:FKPP}
\left\{
\begin{aligned}
    \frac{\partial u}{\partial t}
  - \frac{\partial^2 u}{\partial x^2}
& = f(u) \quad\mbox{ for }\, (x,t)\in \RR\times \RR_+ \,;
\\
  u(x,0)&= u_0(x) \quad\mbox{ for }\, x\in \RR \,.
\end{aligned}
\right.
\end{equation}
This equation was derived by
{\sc R.~A.\ Fisher} \cite{Fisher} in $1937$
and first mathematically analyzed by
{\sc A.\ Kolmogorov}, {\sc I.\ Petrovski}, and {\sc N.\ Piscounov}
\cite{KPP} in the same year.
This is a mathematical problem originating in a morphogenesis model
described in 
{\sc J.~D.\ Murray} \cite{Murray}, {\S}13.3, pp.\ 444--449.
More precisely, we are interested in the convergence of a solution
$u(x,t)$ to a {\it travelling wave\/}
$U(x-ct)$ as $t\to +\infty$, uniformly for $x\in \RR$.
Of course, $c\in \RR$ stands for the {\it speed\/} of the travelling wave.
In particular, we investigate
the {\it stability\/} of travelling waves.
We assume $f(0) = f(1) = 0$ and focus on
biologically (or physically) meaningful solutions satisfying
$0\leq u(x,t)\leq 1$ as $u$ corresponds to the ratio of
an advantageous gene in a population.
Obviously, the ``extreme'' cases $u\equiv 0$ or $u\equiv 1$
throughout $\RR\times \RR_+$ are trivial solutions.
Consequently, we will be interested primarily in nontrivial solutions
$0\leq u(x,t)\leq 1$, with $u\not\equiv 0$ and $u\not\equiv 1$,
and in nontrivial travelling waves
$0\leq U(x-ct)\leq 1$ connecting the two trivial equilibrium states,
i.e.,
\begin{equation}
\label{lim:FKPP}
  \lim_{z\to -\infty} U(z) = 0 \quad\mbox{ and }\quad
  \lim_{z\to +\infty} U(z) = 1 \,.
\end{equation}

The long\--time asymptotic behavior of the dynamical system
generated by the Fisher\--KPP equation \eqref{e:FKPP} is studied, e.g.,
in {\sc P.~C.\ Fife} and {\sc J.~B.\ Mc{L}eod} \cite{Fife-McLeod}
under the standard hypothesis of $f\colon \RR\to \RR$
continuously differentiable (i.e., of class $C^1$).
This crucial hypothesis guarantees not only the existence of
a unique classical solution $u\colon \RR\times \RR_+\to [0,1]$
of the Cauchy (initial value) problem \eqref{e:FKPP},
but facilitates also the linearization of this problem about
a travelling wave $U(x-ct)$.
Such solutions generate a ``smooth'' ($C^1$) dynamical system
whose properties can be investigated by well\--known methods;
see, e.g.,
{\sc J.~K.\ Hale} \cite{Hale} or {\sc D.\ Henry} \cite{Henry}.
Further studies of convergence of solutions to a travelling wave
and front propagation in the semilinear Fisher\--KPP equation with
a $C^1$-reaction function $f$ can be found in
{\sc D.~G.\ Aronson} and {\sc H.~F.\ Weinberger} \cite{AronWein},
{\sc F.\ Hamel} and {\sc N.\ Nadirashvili} \cite{Hamel-Nadira}, and
{\sc H.\ Matano} and {\sc T.\ Ogiwara} \cite{Matano-Og-1, Matano-Og-2}.
Unlike in \cite{Matano-Og-1, Matano-Og-2},
we do not impose any stability condition on a travelling wave
$U\colon \RR\to \RR$, nor do we assume that $f$ is $C^1$,
to obtain its monotonicity.
Thanks to our ``global'' conditions on $f$,
we are able to prove the monotonicity of~$U$
directly from the equation.

In our present work we relax the differentiability hypothesis on $f$
to being only $\alpha$-H\"older\--continuous ($0 < \alpha < 1$)
and ``one\--sided'' Lipschitz\--continuous
(i.e.,
\begin{math}
  s\mapsto f(s) - Ls\colon \RR\to \RR
\end{math}
is monotone decreasing, for some constant $L\in \RR_+$).
In particular, our hypotheses allow for the singular derivatives
\begin{equation}
\label{e:f'(0),f'(1)}
  f'(0) = \lim_{s\to 0} \frac{f(s)}{s} = -\infty
    \quad\mbox{ and }\quad
  f'(1) = \lim_{s\to 1} \frac{f(s)}{s-1} = -\infty \,.
\end{equation}
This type of a reaction function $f$ has been studied extensively
in biological models of various kinds of {\it logistic growth\/}
in {\sc A.\ Tsoularis} and {\sc J.\ Wallace} \cite{Tsoular_Wall}.

The situation sketched in \eqref{e:f'(0),f'(1)}
excludes application of the standard linearization procedure
about a travelling wave widely used in
\cite{Fife-McLeod, Matano-Og-1, Matano-Og-2, Murray}.
Our weaker differentiability conditions on $f$ allow for
a wider and, perhaps, also more realistic class of travelling waves.
In contrast to the usual type of travelling waves $U(x-ct)$
that appear in \cite{Fife-McLeod, Murray} with $U'(z) > 0$
for all $z\in \RR$ and the limits \eqref{lim:FKPP}, we allow also for
\begin{equation}
\label{e:lim_FKPP}
\left\{
\begin{array}{cl}
  U(z) = 0 &\quad\mbox{ if }\, -\infty < z\leq z_0 \,;
\\
  U'(z) > 0 &\quad\mbox{ if }\, z_0 < z < z_1 \,;
\\
  U(z) = 1 &\quad\mbox{ if }\, z_1\leq z < \infty \,,
\end{array}
\right.
\end{equation}
for some $-\infty < z_0 < z_1 < \infty$.
Mathematically, this phenomenon may occur due to $f(s)$
not being Lipschitz\--continuous at the equilibrium points
$s=0$ and $s=1$; cf.\ \eqref{e:f'(0),f'(1)}
and Example~\ref{exam-z_0<z_1}.
Since the linearization procedure is not available,
we are forced to develop alternative methods to treat
the dynamical system generated by \eqref{e:FKPP} and to investigate
its long\--time asymptotic behavior.
Whereas the uniqueness of a solution
$u\colon \RR\times \RR_+\to [0,1]$ to the Cauchy problem \eqref{e:FKPP}
is guaranteed by the hypothesis on
``one\--sided'' Lipschitz continuity of $f$,
the family of all possible travelling waves has to be retrieved by
entirely different methods partially developed already in 
{\sc P. Dr\'abek}, {\sc R.~F. Man\'asevich}, and {\sc P. Tak\'a\v{c}}
\cite{DrabManTak} and
{\sc P. Dr\'abek} and {\sc P. Tak\'a\v{c}} \cite{DrabTak-1}
for $f$ lacking the Lipschitz continuity.
The uniqueness of their speed $c$ is proved in
\cite[Theorem 3.1]{DrabTak-1}, whereas
the monotonicity and uniqueness (up to a spatial translation)
of their profile $U$ will be proved in our present article
(see {\rm Proposition~\ref{prop-uniq_TW}, Part~{\rm (b)}}).

Before stating our main convergence result below,
we would like to remark that
our hypotheses on the reaction function $f\colon \RR\to \RR$
stated in Sections \ref{s:Trav_Waves} and~\ref{s:exist:FKPP},
{\bf (H1)} and {\bf (H2)}, guarantee the uniqueness of the speed
$c\in \RR\setminus \{ 0\}$ and the profile $U\colon \RR\to \RR$
(up to a spatial translation);
see Proposition~\ref{prop-uniq_TW}, Parts {\rm (b)} and~{\rm (d)}.
Adding another hypothesis on $f$ in {\S}\ref{ss:sub-sup},
{\bf (H3)}, we will prove in
Proposition~\ref{prop-weak_sol} that
the Cauchy problem \eqref{e:FKPP} possesses a unique mild solution
$u\colon \RR\times \RR_+\to \RR$
for any initial data $u_0\in L^{\infty}(\RR)$ satisfying
$0\leq u_0\leq 1$ in~$\RR$.
Our main result can be stated loosely as follows:

\begin{theorem}\label{thm-conv_TW}
{\rm (Convergence to a travelling wave.)}$\,$
Let\/ $f\colon \RR\to \RR$ be H\"older\--continuous and\/
$f(0) = f(s_0) = f(1) = 0$ for some $0 < s_0 < 1$.
Assume that\/ $f$ satisfies some additional ``technical'' hypotheses
that guarantee also
$-\infty < z_0 < z_1 < +\infty$.
Finally, assume that the initial data $u_0\colon \RR\to \RR$
are Lebesgue\--measurable and satisfy\hfil\break
$0\leq u_0(x)\leq 1$ for every $x\in \RR$, together with
\begin{equation*}
    \limsup_{x\to -\infty} u_0(x) < s_0
  < \liminf_{x\to +\infty} u_0(x) \,,
\end{equation*}
meaning
\begin{equation*}
    \lim_{n\to \infty}\;
    \esssup_{-\infty < x\leq -n} u_0(x) < s_0
  < \lim_{n\to \infty}\;
    \essinf_{n\leq x < +\infty}  u_0(x) \,.
\end{equation*}
Then the (unique mild) solution
$u\colon \RR\times \RR_+\to \RR$ of problem \eqref{e:FKPP} satisfies
$0\leq u(x,t)\leq 1$ for all\/ $(x,t)\in \RR\times \RR_+$
and there exists a unique spatial shift $\zeta\in \RR$ such that\/
\begin{equation}
\label{e:u(x,t)->U(x)}
  \sup_{x\in \RR}\
  | u(x,t) - U(x-ct + \zeta) | \,\longrightarrow\, 0
  \quad\mbox{ as }\, t\to \infty \,.
\end{equation}
\end{theorem}
\par\vskip 10pt

This theorem is easily derived from Theorem~\ref{thm-Lyapunov}
proved in Section~\ref{s:Time-Asympt}
using the transformation with the moving coordinate $z = x-ct$.

This article is organized as follows.
In the next section (Section~\ref{s:Trav_Waves})
we analyze the travelling waves.
The existence and uniqueness of a solution to
the Cauchy (initial value) problem \eqref{e:FKPP}
are established in Section~\ref{s:exist:FKPP}.
In this section we prove also
the usual weak comparison principle for parabolic problems
({\S}\ref{ss:sub-sup}).
Moreover, we construct a special pair of ordered sub- and super\-solutions
to the Cauchy problem \eqref{e:FKPP}
by modifying the travelling waves ({\S}\ref{ss:Trav-sub/super}).
This pair is subsequently used in Section~\ref{s:Time-Asympt}
to prove the Lyapunov stability of a travelling wave
({\S}\ref{ss:Stability})
as well as the convergence of a solution $u(x,t)$
to a single travelling wave $U(x-ct + \zeta)$ as $t\to \infty$,
uniformly for $x\in \RR$; cf.\ Theorem~\ref{thm-Lyapunov}.

\section{Analysis of Travelling Waves}
\label{s:Trav_Waves}

We reformulate the Cauchy problem \eqref{e:FKPP} for $u(x,t)$
as an equivalent initial value problem for the unknown function
$v(z,t) = u(z+ct,t)\equiv u(x,t)$ with the moving coordinate $z = x-ct$:
\begin{equation}
\label{e_c:FKPP}
\left\{
\begin{aligned}
    \frac{\partial v}{\partial t}
  - \frac{\partial^2 v}{\partial z^2}
  - c\, \frac{\partial v}{\partial z}
& = f(v) \,,
    \quad (z,t)\in \RR\times \RR_+ \,;
\\
  v(z,0)&= v_0(z) \,.
\end{aligned}
\right.
\end{equation}
We will show in Proposition~\ref{prop-monot_TW}
({\S}\ref{ss:U'>0} below)
that every travelling wave $u(x,t) = U(x-ct)$
for \eqref{e:FKPP} must have a monotone increasing profile
$U\colon \RR\to \RR$ satisfying \eqref{lim:FKPP}
and $U'(z) > 0$ for every $z\in \RR$
such that $0 < U(z) < 1$.
Hence, by \cite[Theorem 3.1]{DrabTak-1},
the speed $c$ of a travelling wave $U(x-ct)$ in eq.~\eqref{e:FKPP}
is determined uniquely by~$f$.
This means that we are able to ``reduce'' the problem of finding
all {\em travelling waves\/} $u(x,t) = U(x-ct)$ for \eqref{e:FKPP}
to the problem of investigating all {\em stationary solutions\/}
$v(z,t) = U(z)$ for \eqref{e_c:FKPP} that obey also \eqref{lim:FKPP}.
Of course, $v_0\equiv u_0$ in~$\RR$.
Since both these families are described by the same profile
$U\colon \RR\to [0,1]$ that obeys \eqref{lim:FKPP},
the uniqueness of which (up to a spatial translation)
will be proved later in
{\rm Proposition~\ref{prop-uniq_TW}, Part~{\rm (b)}},
we introduce the following terminology for the latter:
We call $v(z,t) = U(z)$ and any other stationary solution
$v(z,t) = U(z + \zeta)$, $\zeta\in \RR$ -- a constant,
that obeys \eqref{lim:FKPP},
a {\em\bfseries TW\--solution\/} of problem \eqref{e_c:FKPP}.
It should not to be confused with the travelling wave
$u(x,t) = U(x-ct + \zeta)$ in the original problem \eqref{e:FKPP}.
In fact, to identify all possible TW\--solutions of problem
\eqref{e_c:FKPP}, \eqref{lim:FKPP},
we first establish their monotonicity
(in Proposition~\ref{prop-monot_TW}),
then prove the uniqueness (up to a spatial translation)
of their profile $U\colon \RR\to [0,1]$ satisfying \eqref{lim:FKPP}
(in {\rm Proposition~\ref{prop-uniq_TW}, Part~{\rm (b)}}).
Finally, all possible TW\--solutions are described in
{\rm Proposition~\ref{prop-uniq_TW}, Part~{\rm (d)}}.
Thus, except when we wish to stress
the biological (or physical) meaning of travelling waves,
we prefer to treat the Cauchy problem \eqref{e_c:FKPP}
and investigate its TW\--solutions, rather than problem \eqref{e:FKPP}.

As usual, we denote
$\RR\eqdef (-\infty, \infty)$, $\RR_+\eqdef [0,\infty)$,
and assume that the reaction term $f$ satisfies
the following basic hypotheses:
\begin{enumerate}
\renewcommand{\labelenumi}{{\bf (H\alph{enumi})}}
\item[{\bf (H1)}]
$f\colon \RR\to \RR$ is a continuous, but
{\it not necessarily smooth\/} function,
such that
$f(0) = f(s_0) = f(1) = 0$ for some $0 < s_0 < 1$, together with
$f(s) < 0$ for every $s\in (0,s_0)\cup (1,\infty)$,
$f(s) > 0$ for every $s\in (-\infty,0)\cup (s_0,1)$, and
\begin{equation}
\label{int:f(u)}
  F(r)\eqdef \int_0^r f(s) \,\mathrm{d}s < 0
    \quad\mbox{ whenever }\, 0 < r\leq 1 \,.
\end{equation}
\end{enumerate}
We remark that the case $F(1) = 0$ would prevent the existence of
any travelling wave with speed $c\neq 0$; cf.\ eq.~\eqref{e:id_c}.
That is why we assume $F(1) < 0$.
We will see in our proof of Proposition~\ref{prop-uniq_TW}
(in {\S}\ref{ss:Trav-ODE} below)
that {\rm Hypothesis {\bf (H1)}} is sufficient to determine
the speed $c$ and the profile $U\colon \RR\to [0,1]$ uniquely,
the latter up to a spatial translation.
The reaction function $f$ satisfying {\bf (H1)} models the so\--called
{\em heterozygote inferior\/} case in the genetic model studied, e.g.,
in {\sc D.~G.\ Aronson} and {\sc H.~F.\ Weinberger}
\cite[p.~35]{AronWein}.

In the next paragraph, {\S}\ref{ss:U'>0},
we will show that
$U\colon \RR\to \RR$ must be monotone increasing with
$U'> 0$ on a suitable open interval $(z_0,z_1)\subset \RR$, such that
\begin{equation}
\label{z_0,1:FKPP}
  \lim_{z\to z_0+} U(z) = 0 \quad\mbox{ and }\quad
  \lim_{z\to z_1-} U(z) = 1 \,.
\end{equation}
We would like to remark that the cases of
$z_0 > -\infty$ and/or $z_1 < +\infty$
render qualitatively different travelling waves than
the classical case $(z_0,z_1) = \RR$ which has been studied
in the original works \cite{Fisher, KPP} and in the literature
(\cite{AronWein, Fife-McLeod, Hamel-Nadira, Murray, Murray-I}).
Using this setting, in paragraph {\S}\ref{ss:Trav-ODE}
we are able to find a {\it first integral\/} for
the second\--order equation for $U$:
\begin{equation}
\label{eq:FKPP}
    \frac{\mathrm{d}^2 U}{\mathrm{d}z^2}
  + c\, \frac{\mathrm{d}U}{\mathrm{d}z} + f(U)
  = 0 \,,
    \quad z\in (z_0,z_1) \,.
\end{equation}
%

\subsection{Monotonicity of a Travelling Wave}
\label{ss:U'>0}

In accordance with problem \eqref{eq:FKPP}, \eqref{lim:FKPP},
we investigate the following problem for an unknown $C^2$-function
$v\colon \RR\to \RR$:
\begin{align}
\label{eq:v(z)}
  v''(z) + c\, v'(z) + f(v(z)) = 0 \,,\quad z\in \RR \,;
\\
\label{lim:v(z)}
  \lim_{z\to -\infty} v(z) = 0 \quad\mbox{ and }\quad
  \lim_{z\to +\infty} v(z) = 1 \,.
\end{align}
We would like to stress that here we do not impose any hypothesis on
the reaction function $f\colon \RR\to \RR$
that would require some kind of H\"older or
Lipschitz continuity, not even one\--sided.
Thus, we cannot use standard ``local'' uniqueness arguments for
eq.~\eqref{eq:v(z)}
to find the family of all $C^2$-solutions $v\colon \RR\to \RR$
to problem \eqref{eq:v(z)}, \eqref{lim:v(z)}.
Nevertheless, we will be able to demonstrate that,
if the constant $c\in \RR$ is such that at least one solution $U$
to this problem exists, then any other solution
$v\colon \RR\to \RR$
can be obtained by a translation in the argument of $U$
by a fixed number $\zeta\in \RR$, i.e.,
$v(z) = U(z + \zeta)$ for $z\in \RR$; cf.\
Proposition~\ref{prop-uniq_TW} ({\S}\ref{ss:Trav-ODE}) below.
In order to treat questions related to the uniqueness of a solution $v$,
we take advantage of the global properties of $v$,
such as the (strict) monotonicity established in the next proposition,
that are consequences of the global behavior of the reaction function
$f\colon \RR\to \RR$ specified in {\rm Hypothesis {\bf (H1)}}.

Before beginning to investigate rather subtle properties of
travelling waves, we recall the following well\--known identity
which holds for any solution $v\colon \RR\to \RR$
to problem \eqref{eq:v(z)}, \eqref{lim:v(z)}:
\begin{equation}
\label{e:id_c}
    c\int_{-\infty}^{+\infty} v'(z)^2 \,\mathrm{d}z
  = F(0) - F(1) = - F(1)\quad ({} > 0) \,.
\end{equation}
Consequently, $c > 0$.
This identity is obtained by calculating the first integral for
\eqref{eq:v(z)},
\begin{equation}
\label{eq:id_c}
    \frac{1}{2}\, \frac{\mathrm{d}}{\mathrm{d}z} (v'(z)^2)
  + c\, v'(z)^2 + \frac{\mathrm{d}}{\mathrm{d}z} F(v(z)) = 0 \,,\quad
    z\in \RR \,,
\end{equation}
integrating it over a suitable sequence of intervals
$(x_k,y_k)\subset \RR$; $k=1,2,3,\dots$, with
$- x_k, y_k\in (k,k+1)$,
\begin{equation}
\label{int:id_c}
    \frac{1}{2}\, (v'(y_k)^2) - \frac{1}{2}\, (v'(x_k)^2)
  + c\int_{x_k}^{y_k} v'(z)^2 \,\mathrm{d}z
  + F(v(y_k)) - F(v(x_k)) = 0 \,,
\end{equation}
and letting $k\to \infty$.
Here, by the mean value theorem, we construct
$x_k\in (-k-1,-k)$ and $y_k\in (k,k+1)$ such that
$v'(x_k) = v(-k) - v(-k-1)$ and $v'(y_k) = v(k+1) - v(k)$,
respectively.
Thanks to the limits in \eqref{lim:v(z)}, we arrive at
$\lim_{k\to \infty} v'(x_k) = \lim_{k\to \infty} v'(y_k)$ $= 0$.
Hence, \eqref{e:id_c} follows from
\eqref{lim:v(z)} and \eqref{int:id_c}, by letting $k\to \infty$.

We have the following result.

\begin{proposition}\label{prop-monot_TW}
{\rm (Monotonicity.)}$\,$
Let\/ $c\in (0,\infty)$ and assume that\/
$f\colon \RR\to \RR$ satisfies {\rm Hypothesis {\bf (H1)}}.
Then every\/ $C^2$-solution $v\colon \RR\to \RR$
of eq.~\eqref{eq:v(z)} with the limits \eqref{lim:v(z)} satisfies
$0\leq v(z)\leq 1$ and\/ $v'(z)\geq 0$ for all\/ $z\in \RR$.
Moreover, there is an open interval\/ $(z_0,z_1)\subset \RR$,
$-\infty\leq z_0 < z_1\leq +\infty$, such that\/
$v'> 0$ on $(z_0,z_1)$ together with
\begin{equation}
\label{e:z<0,z>1}
\left\{\quad
  \begin{alignedat}{2}
& \lim_{z\to z_0+} v(z) = 0 \quad\mbox{ and }\quad
&&  v(z) = 0 \;\mbox{ if }\; -\infty < z\leq z_0 \,,
\\
& \lim_{z\to z_1-} v(z) = 1 \quad\mbox{ and }\quad
&&  v(z) = 1 \;\mbox{ if }\; z_1\leq z < +\infty \,.
  \end{alignedat}
\right.
\end{equation}
\end{proposition}

\par\vskip 10pt
\proof
We verify our first claim,
$0\leq v(z)\leq 1$ for all $z\in \RR$, by contradiction.
Suppose there is a number $\xi\in \RR$ such that $v(\xi) < 0$.
We make use of the limits in~\eqref{lim:v(z)}
to conclude that there are numbers $\xi_1, \xi_2\in \RR$ such that
$\xi_1 < \xi < \xi_2$ and
$v(\xi) < \min\{ v(\xi_1), v(\xi_2)\}$.
Denoting by $\xi_0\in [\xi_1, \xi_2]$
the (global) minimizer for the function $v$ over the compact interval
$[\xi_1, \xi_2]$, we arrive at $\xi_0\in (\xi_1, \xi_2)$,
$v(\xi_0)\leq v(\xi) < 0$, $v'(\xi_0) = 0$, and $v''(\xi_0)\geq 0$.
Hence, we have also $f(v(\xi_0)) > 0$,
by {\rm Hypothesis {\bf (H1)}}.
We insert these properties of $v(z)$ at $z = \xi_0$
into the left\--hand side of eq.~\eqref{eq:v(z)}, thus arriving at
\begin{equation*}
  v''(\xi_0) + c\, v'(\xi_0) + f(v(\xi_0)) > 0 \,,
\end{equation*}
a contradiction to eq.~\eqref{eq:v(z)}.

We proceed analogously if $\xi\in \RR$ is such that $v(\xi) > 1$.
The limits in~\eqref{lim:v(z)} guarantee that
there are numbers $\xi_1, \xi_2\in \RR$ such that
$\xi_1 < \xi < \xi_2$ and
$v(\xi) > \max\{ v(\xi_1), v(\xi_2)\}$.
Denoting by $\xi_0\in [\xi_1, \xi_2]$
the maximizer for the function $v$ over the compact interval
$[\xi_1, \xi_2]$, we arrive at $\xi_0\in (\xi_1, \xi_2)$,
$v(\xi_0)\geq v(\xi) > 1$, $v'(\xi_0) = 0$, and $v''(\xi_0)\leq 0$.
Hence, we have also $f(v(\xi_0)) < 0$,
by {\rm Hypothesis {\bf (H1)}}.
Similarly as above, we arrive at
\begin{equation*}
  v''(\xi_0) + c\, v'(\xi_0) + f(v(\xi_0)) < 0 \,,
\end{equation*}
a contradiction to eq.~\eqref{eq:v(z)} again.
We have verified that $0\leq v(z)\leq 1$ holds for all $z\in \RR$.

Now it follows from \eqref{lim:v(z)} again that
there is a maximal (nonempty) open interval
$(z_0,z_1)\subset \RR$ such that
$v'(z) > 0$ holds for all $z\in (z_0,z_1)$.
The strict monotonicity of $v$ in $(z_0,z_1)$ guarantees
the existence of the (monotone) limits
\begin{equation}
\label{lim:l_0<v<l_1}
  0\leq \ell_0\eqdef \lim_{z\to z_0+} v(z)
      < \ell_1\eqdef \lim_{z\to z_1-} v(z)\leq 1 \,.
\end{equation}
The maximality of the interval $(z_0,z_1)$ forces either
$(z_0,z_1) = \RR$ in which case the proposition is proved, or else
at least one of the inequalities
$z_0 > -\infty$ and $z_1 < +\infty$ must be valid:
If $z_0 > -\infty$ then $v'(z_0) = 0$, whereas
if $z_1 < +\infty$ then $v'(z_1) = 0$.

We begin to investigate these two possibilities for
the following two extremal cases:
\begin{itemize}
\item[{\rm (i)}]
$\;$
$z_0 > -\infty$, $v'(z_0) = 0$, and $v(z_0) = 0$;
\item[{\rm (ii)}]
$\;$
$z_1 < +\infty$, $v'(z_1) = 0$, and $v(z_1) = 1$.
\end{itemize}

To treat these cases, we introduce the Lyapunov function
$V\equiv V[v]\colon \RR\to \RR$ by
\begin{equation}
\label{e:Lyapunov}
  V(z)\equiv V[v](z)\eqdef
  \frac{1}{2}\, [v'(z)]^2 + F(v(z)) \,,\quad z\in \RR \,.
\end{equation}
Taking advantage of eq.~\eqref{eq:v(z)}, we calculate
\begin{equation}
\label{ineq:Lyapunov}
    \frac{\mathrm{d}}{\mathrm{d}z}\, V(z)
  = v'\, v'' + f(v(z))\, v' = {}- c\, [v'(z)]^2\leq 0
    \quad\mbox{ for all }\, z\in \RR \,.
\end{equation}
In particular, if $V(\zeta_1)\leq V(\zeta_2)$ for some
$-\infty < \zeta_1 < \zeta_2 < +\infty$, then we have
$v'(z)\equiv 0$ and
$v(z)\equiv v(\zeta_1) = v(\zeta_2)$ for every
$z\in [\zeta_1, \zeta_2]$.

{\it Case}~{\rm (i)}.$\;$
We claim that $v(z_0) = 0$ implies $v(z) = 0$ for all $z\leq z_0$.
$\;${\it Proof:\/}
Since $0\leq v(z)\leq 1$ holds for all $z\in \RR$, we have also
$v'(z_0) = 0$.
Thus, on the contrary to our claim, suppose that
$v(\hat{z}) > 0$ at some $\hat{z}\in (-\infty,z_0)$.
Taking into account also
$\lim_{z\to -\infty} v(z) = 0$, we conclude that
the function $v$ attains its (global) maximum
$v(\xi)$ over the interval $(-\infty,z_0]$ at some point
$\xi\in (-\infty,z_0)$.
Consequently, we have $v'(\xi) = 0$ and
$v(z_0) = 0 < v(\hat{z})\leq v(\xi)\leq 1$.
We make use of standard properties of the Lyapunov function $V$
to compare
\begin{equation*}
  F(v(\xi)) = V(\xi) > V(z_0) = F(v(z_0)) = F(0) = 0
\end{equation*}
which contradicts our hypothesis \eqref{int:f(u)}, i.e.,
$F(r) < 0$ with $r = v(\xi)\in (0,1]$.
We have proved that $v(z) = 0$ holds for all $z\leq z_0$.

{\it Case}~{\rm (ii)}.$\;$
We claim that $v(z_1) = 1$ implies $v(z) = 1$ for all $z\geq z_1$.
$\;${\it Proof:\/}
Again, from $0\leq v\leq 1$ we deduce also $v'(z_1) = 0$.
Thus, on the contrary to our claim, suppose that
$v(\hat{z}) < 1$ at some $\hat{z}\in (z_1,+\infty)$.
Then standard properties of the Lyapunov function $V$ yield
\begin{equation*}
  0 > F(1) = F(v(z_1)) = V(z_1) > V(\hat{z})\geq V(z)\geq F(v(z))
    \quad\mbox{ for all }\, z\geq \hat{z} \,.
\end{equation*}
Letting $z\to +\infty$ and using
$\lim_{z\to +\infty} v(z) = 1$, we arrive at
\begin{equation*}
  F(1) > V(\hat{z})\geq \lim_{z\to +\infty} F(v(z)) = F(1) \,,
\end{equation*}
a contradiction again.
We have verified $v(z) = 1$ for all $z\geq z_1$.

In our next step we show that there is no number $\xi\in \RR$ such that
$v'(\xi) = 0$ and $0 < v(\xi) < 1$.
On the contrary, if $\xi\in \RR$ is such that
$v'(\xi) = 0$ and $0 < v(\xi) < 1$,
then we distinguish between the following two alternatives:

{\it Alt.~$1$.}$\;$ Suppose that
$(0 <)$ $s_0\leq v(\xi) < 1$ and $v'(\xi) = 0$.
By {\rm Hypothesis {\bf (H1)}}, the graph of the function
\begin{equation*}
  F(r)\eqdef \int_0^r f(s) \,\mathrm{d}s
    \quad\mbox{ for }\, r\in \RR
\end{equation*}
shows that there is a unique number $s_1\in (0,s_0)$ such that
$F(s_1) = F(v(\xi))$, thanks to
$F(s_0)\leq F(v(\xi)) < F(1) < 0 = F(0)$.
From
\begin{equation*}
  F(v(\xi)) = V(\xi)\geq V(z)\geq F(v(z))
    \quad\mbox{ for all }\, z\geq \xi
\end{equation*}
and the strict monotonicity of the function $F$ on the intervals
$[s_1,s_0]$ ($F$ strictly monotone decreasing) and
$[s_0,1]$ ($F$ strictly monotone increasing) we deduce that
\begin{equation*}
  0 < s_1\leq v(z)\leq v(\xi) < 1
    \quad\mbox{ for all }\, z\geq \xi \,.
\end{equation*}
But this is a contradiction to
$\lim_{z\to +\infty} v(z) = 1$.

{\it Alt.~$2$.}$\;$ Now suppose that
$0 <  v(\xi) < s_0$ $(< 1)$ and $v'(\xi) = 0$.
Eq.~\eqref{eq:v(z)} being equivalent with
\begin{equation}
\label{equiv:v(z)}
    \ee^{-cz}\cdot \frac{\mathrm{d}}{\mathrm{d}z}
    \left( \ee^{cz}\, v'(z)\right)
  = {}- f(v(z)) \,,\quad z\in \RR \,,
\end{equation}
we conclude that the function
\begin{math}
  z\mapsto \ee^{cz}\, v'(z) \colon \RR\to \RR
\end{math}
is strictly monotone increasing on the maximal open interval
$J\subset \RR$ such that
\begin{equation*}
  \xi\in J\subset
  v^{-1}(0,s_0)\eqdef \{ z\in \RR\colon 0 < v(z) < s_0\} \,.
\end{equation*}
Consequently,
\begin{equation*}
  \ee^{cz}\, v'(z) < \ee^{c\xi}\, v'(\xi) = 0
    \quad\mbox{ holds for all }\, z\in J\cap (-\infty, \xi) \,,
\end{equation*}
which gives
$J\cap (-\infty, \xi]\cap (z_0,z_1) = \emptyset$ and
$0 < v(\xi) < v(z) < s_0$ for all $z\in J\cap (-\infty, \xi)$.
Similarly,
\begin{equation*}
  \ee^{cz}\, v'(z) > \ee^{c\xi}\, v'(\xi) = 0
    \quad\mbox{ for all }\, z\in J\cap (\xi, +\infty) \,,
\end{equation*}
which gives
$0 < v(\xi) < v(z) < s_0$ for all $z\in J\cap (\xi, +\infty)$.
Hence, $\xi\in J$ is a strict (global) minimizer for the function $v$
over the interval $J$.
As a consequence,
$0 < v(\xi)\leq v(z) < s_0 < 1$ holds for all $z\in J$.
The limits in~\eqref{lim:v(z)} force
$-\infty < \inf J < \sup J < +\infty$, so that
$J = (z_0',z_1')$ is a bounded open interval and
$v(z_0') = v(z_1') = s_0$ together with
$v'(z_0') < v'(\xi) = 0 < v'(z_1')$.
Recalling
$\lim_{z\to -\infty} v(z) = 0$, we now observe that
$v$ attains its maximum over the interval
$(-\infty, z_0')$ at some $\xi'\in (-\infty, z_0')$, i.e.,
$0\leq v(z)\leq v(\xi')$ for all $z\leq z_0'$, together with
$s_0 < v(\xi')\leq 1$ and $v'(\xi') = 0$.

By {\rm Case~(ii)} above, the alternative $v(\xi') = 1$ would force
$v(z)\equiv 1$ for all $z\geq \xi'$.
We may set $z = z_0'> \xi'$ to arrive at the following contradiction,
\begin{math}
  v(z_0') = 1 > s_0 = v(z_0') .
\end{math}
Consequently, we must have $s_0 < v(\xi') < 1$ and $v'(\xi') = 0$.
But this situation for $\xi'\in \RR$ corresponds to
Alternative~$1$ for $\xi\in \RR$ above which has already been excluded.

We combine Alternatives $1$ and~$2$ to conclude that
if $\xi\in \RR$ is a number with $v'(\xi) = 0$ then we have
either $v(\xi) = 0$ or else $v(\xi) = 1$.
Furthermore, by Cases {\rm (i)} and~{\rm (ii)}
(treated before Alternatives $1$ and~$2$), we have even
$v(z)\equiv v(\xi) = 0$ for all $z\leq \xi$ (in Case~{\rm (i)}), or
$v(z)\equiv v(\xi) = 1$ for all $z\geq \xi$ (in Case~{\rm (ii)}).

Now we are ready to complete our proof.
It remains to treat the following two possibilities:
{\rm (a)}$\;$
$z_0 > -\infty$ and $v'(z_0) = 0$, and
{\rm (b)}$\;$
$z_1 < +\infty$ and $v'(z_1) = 0$.

We have
$0\leq v(z_0) < v(z) < v(z_1)\leq 1$ for all $z\in (z_0,z_1)$.
We can have neither $0 < v(z_0) < 1$ nor $0 < v(z_1) < 1$,
by Alternatives $1$ and~$2$.

Concerning {\rm (a)}, we must have $v(z_0) = 0$ and, thus,
$v(z)\equiv v(z_0) = 0$ for all $z\leq z_0$, as well, by Case~{\rm (i)}.

Similarly, concerning {\rm (b)}, we must have $v(z_1) = 1$ and, thus,
$v(z)\equiv v(z_1) = 1$ for all $z\geq z_1$, again, by Case~{\rm (ii)}.

To conclude our proof, we recall that
$v'(z) > 0$ holds for every $z\in (z_0,z_1)$.
\qed
\par\vskip 10pt

\begin{remark}\label{rem-monot_TW}\nopagebreak
\begingroup\rm
Let $v\colon \RR\to \RR$ and $(z_0,z_1)\subset \RR$
be as in {\rm Proposition~\ref{prop-monot_TW}} above.
Since our problem \eqref{eq:v(z)}, \eqref{lim:v(z)}
is invariant with respect to a shift
$z\mapsto z + \zeta$ in the variable $z\in \RR$,
by a fixed number $\zeta\in \RR$, then also the function
$\tilde{v}\colon \RR\to \RR$,
$\tilde{v}(z)\eqdef v(z + \zeta)$ for $z\in \RR$,
is a $C^2$-solution of eq.~\eqref{eq:v(z)} satisfying
$0\leq \tilde{v}\leq 1$ and $\tilde{v}'\geq 0$ on $\RR$,
$\tilde{v}'> 0$ on
$(\tilde{z}_0, \tilde{z}_1) = (z_0 - \zeta, z_1 - \zeta)$, and
\begin{equation}
\label{e_zeta:z<0,z>1}
\left\{\quad
  \begin{alignedat}{2}
& \lim_{z\to \tilde{z}_0+} \tilde{v}(z) = 0 \quad\mbox{ and }\quad
&&  \tilde{v}(z) = 0 \;\mbox{ if }\; -\infty < z\leq \tilde{z}_0 \,,
\\
& \lim_{z\to \tilde{z}_1-} \tilde{v}(z) = 1 \quad\mbox{ and }\quad
&&  \tilde{v}(z) = 1 \;\mbox{ if }\; \tilde{z}_1\leq z < +\infty \,.
  \end{alignedat}
\right.
\end{equation}
We will see in the next paragraph
({\S}\ref{ss:Trav-ODE}, Proposition~\ref{prop-uniq_TW})
that problem \eqref{eq:v(z)}, \eqref{lim:v(z)}
has a $C^2$-solution $v\colon \RR\to \RR$
for precisely one value of the speed $c\in \RR$.
Furthermore, this solution is monotone increasing and
unique up to a shift by a constant $\zeta\in \RR$, as described above.
\endgroup
\end{remark}

\subsection{Reduction to a First\--Order O.D.E.}
\label{ss:Trav-ODE}

We use a phase plane transformation
(cf.\ {\sc J.~D.\ Murray} \cite{Murray-I}, {\S}13.2, pp.\ 440--441,
 and {\sc P.\ Dr\'abek} and {\sc P.\ Tak\'a\v{c}} \cite{DrabTak-1})
in order to describe all monotone increasing travelling waves
$u(x,t) = v(x-ct,t)\equiv U(x-ct - \zeta)$
where $U\colon \RR\to \RR$ is a $C^2$-solution of
problem \eqref{eq:v(z)}, \eqref{lim:v(z)}
as specified in {\rm Remark~\ref{rem-monot_TW}}
and normalized by $U(0) = s_0$, and
$\zeta\in \RR$ is a suitable constant.
Instead of trying to find such a unique function
$U = U(z)$ of $z\in \RR$, we will calculate the derivative
${\mathrm{d}z} / {\mathrm{d}U}$ of its inverse function
$U\mapsto z = z(U)$ as a function of $U\in (0,1)$.
In fact, below we find a nonlinear differential equation
for the derivative
\begin{equation*}
  U'(z) = \genfrac{(}{)}{}0{\mathrm{d}z}{\mathrm{d}U}^{-1}
        \equiv \frac{1}{z'(U)}
  \quad\mbox{ as a function of }\, U\in (0,1) \,.
\end{equation*}
To this end, we make the substitution
\begin{equation}
\label{e:dU/dz}
  V\eqdef \genfrac{}{}{}0{\mathrm{d}U}{\mathrm{d}z} > 0
  \quad\mbox{ for }\, z\in (z_0,z_1)
\end{equation}
and consequently look for $V = V(U)$ as a function of
$U\in (0,1)$ that satisfies the following differential equation
obtained from eq.~\eqref{eq:FKPP}:
\begin{equation*}
    \frac{\mathrm{d}V}{\mathrm{d}U}\cdot \frac{\mathrm{d}U}{\mathrm{d}z}
  + c\, \frac{\mathrm{d}U}{\mathrm{d}z} + f(U)
  = 0 \,,
    \quad z\in (z_0,z_1) \,,
\end{equation*}
that is,
\begin{equation}
\label{eq:FKPP:V(U)}
    \frac{\mathrm{d}V}{\mathrm{d}U}\cdot V + c\, V + f(U) = 0 \,,
    \quad U\in (0,1) \,.
\end{equation}
Hence, we are looking for the inverse function
$U\mapsto z(U)$ with the derivative
\begin{equation*}
  \genfrac{}{}{}0{\mathrm{d}z}{\mathrm{d}U} = \frac{1}{V(U)} > 0
  \quad\mbox{ for }\, U\in (0,1) \,, \quad\mbox{ such that }\,
  z(s_0) = 0 \,.
\end{equation*}
Finally, in equation \eqref{eq:FKPP:V(U)} we make the substitution
\begin{equation}
\label{e:y=V^p'}
  y = V^2 = \genfrac{|}{|}{}0{\mathrm{d}U}{\mathrm{d}z}^2 > 0
\end{equation}
and write $r$ in place of $U$, thus arriving at
\begin{equation*}
    \frac{1}{2}\cdot \frac{\mathrm{d}y}{\mathrm{d}r}
  + c\, \sqrt{y} + f(r) = 0 \,,
    \quad r\in (0,1) \,.
\end{equation*}
This means that the unknown function
$y\colon (0,1)\to (0,\infty)$ of $r$ verifies also
\begin{equation}
\label{eq:FKPP:y(r)}
    \frac{\mathrm{d}y}{\mathrm{d}r}
  = {}- 2\left( c\, \sqrt{y^{+}} + f(r) \right) \,,
    \quad r\in (0,1) \,,
\end{equation}
where $y^{+} = \max\{ y,\, 0\}$.
Since we require that $U = U(z)$ be sufficiently smooth,
at least continuously differentiable, with
$U'(z)\to 0$ as $z\to z_0+$ and $z\to z_1-$, the function
$y = y(r) = | \mathrm{d}U / \mathrm{d}z |^2$
must satisfy the boundary conditions
\begin{equation}
\label{bc:FKPP:y(r)}
  y(0) = y(1) = 0 \,.
\end{equation}

We take advantage of Hypotheses {\bf (H1)} on $f$ to claim that,
by a result in
{\sc P.\ Dr\'abek} and {\sc P.\ Tak\'a\v{c}}
\cite[Corollary 5.5]{DrabTak-1},
the overdetermined first\--order, two\--point boundary value problem
\eqref{eq:FKPP:y(r)}, \eqref{bc:FKPP:y(r)}
admits a $C^1$-solution $y\colon [0,1]\to \RR$
for precisely one value of $c\in \RR$;
this solution is positive in $(0,1)$.
We emphasize again that the nonlinearity
$y\mapsto \sqrt{y^{+}}$ in the differential equation
\eqref{eq:FKPP:y(r)} does not satisfy a local Lipschitz condition,
so, due to the lack of uniqueness of a solution,
the standard shooting method cannot be applied directly.
Finally, the desired TW\--solution $U = U(z)$ is calculated from
eq.~\eqref{e:dU/dz} with $V = \sqrt{y}$.
More precisely, the function $U$, which is monotone increasing,
is calculated by integrating the differential equation
\begin{equation}
\label{e:dz/dU}
  \mathrm{d}z =
  \genfrac{}{}{}0{\mathrm{d}U}{\sqrt{y(U)}}
  \quad\mbox{ for }\, U\in (0,1) \,, \quad\mbox{ with }\,
  z(s_0) = 0 \,.
\end{equation}

We summarize our results from paragraphs
{\S}\ref{ss:U'>0} and {\S}\ref{ss:Trav-ODE}
in the following proposition:

\begin{proposition}\label{prop-uniq_TW}
{\rm (TW\--solutions.)}$\,$
Assume that\/
$f\colon \RR\to \RR$ satisfies {\rm Hypothesis {\bf (H1)}}.
Then the following statements are valid:
\begin{itemize}
\item[{\rm (a)}]
The two\--point boundary value problem
\eqref{eq:FKPP:y(r)}, \eqref{bc:FKPP:y(r)}
has a $C^1$-solution $y\colon [0,1]\to \RR$
for precisely one value of the speed $c\in \RR$;
this critical value, denoted by $c^{\ast}$, is positive, $c^{\ast} > 0$.
The corresponding solution, $y\equiv y_{c^{\ast}}$,
is unique and positive in $(0,1)$.
\item[{\rm (b)}]
Problem \eqref{eq:v(z)}, \eqref{lim:v(z)} with $c\in \RR$
has a $C^2$-solution $v\colon \RR\to \RR$
if and only if $c = c^{\ast}$.
This solution, $v\equiv v_{c^{\ast}}$, is unique if we require also
$v(0) = s_0$; we denote it by $U\equiv v_{c^{\ast}}$.
\item[{\rm (c)}]
The solution $U\colon \RR\to \RR$ satisfies
$0\leq U(z)\leq 1$ and $U'(z)\geq 0$ for every\/ $z\in \RR$.
Moreover, there is an open interval\/
$(z_0,z_1)\subset \RR$, $-\infty\leq z_0 < z_1\leq +\infty$,
such that\/
$0 < U(z) < 1$ and $U'(z) > 0$ hold for every\/ $z\in (z_0,z_1)$, and\/
\begin{equation}
\label{e:U:z<0,z>1}
\left\{\quad
  \begin{alignedat}{2}
& \lim_{z\to z_0+} U(z) = 0 \quad\mbox{ and }\quad
&&  U(z) = 0 \;\mbox{ if }\; -\infty < z\leq z_0 \,,
\\
& \lim_{z\to z_1-} U(z) = 1 \quad\mbox{ and }\quad
&&  U(z) = 1 \;\mbox{ if }\; z_1\leq z < +\infty \,.
  \end{alignedat}
\right.
\end{equation}
\item[{\rm (d)}]
Without the condition $v(0) = s_0$, all other solutions
$v\colon \RR\to \RR$
of problem \eqref{eq:v(z)}, \eqref{lim:v(z)} with $c = c^{\ast}$
are given by
$v(z) = U(z + \zeta)$ for $z\in \RR$, where
$\zeta\in \RR$ is arbitrary.
\end{itemize}
\end{proposition}
\par\vskip 10pt

\subsection{The Asymptotic Shape of Travelling Waves}
\label{ss:Trav-Shape}

The asymptotic shape of the TW\--solution $U(z)$ as $z\to \mp\infty$
is determined by the asymptotic behavior of the integral
\begin{equation}
\label{int:dz/dU}
  \int_{r_0}^{r_1}
  \genfrac{}{}{}0{\mathrm{d}r}{V(r)}
    \quad\mbox{ as $r_0\to 0+$ and $r_1\to 1-$, respectively; }
\end{equation}
cf.\ eq.~\eqref{e:dz/dU}.
Indeed, we observe that the integral
\begin{equation}
\label{def:dz/dU}
  z(U) = z(s_0) + \int_{s_0}^{U}
  \genfrac{}{}{}0{\mathrm{d}r}{V(r)} \,,
    \quad\mbox{ for }\, U\in (0,1) \,,
\end{equation}
renders the inverse function of a TW\--solution $U(z)$
that is determined uniquely by the point
$z(s_0) = 0$ at which $U(0) = s_0$.
Next, let us consider the limits
\begin{equation}
\label{lim:dz/dU}
  z_0\eqdef \lim_{U\to 0+} z(U)\geq -\infty \qquad\mbox{ and }\qquad
  z_1\eqdef \lim_{U\to 1-} z(U)\leq +\infty \,.
\end{equation}
These limits can be finite or infinite,
$-\infty\leq z_0 < 0 < z_1\leq +\infty$, depending on whether
the integral in \eqref{int:dz/dU} is convergent or divergent
as $r_0\to 0+$ and $r_1\to 1-$, respectively.
This, in turn, depends on the asymptotic behavior of
the function $f(r)$ as $r\to 0+$ and $r\to 1-$,
thanks to $y = V^2$ being a solution to problem
\eqref{eq:FKPP:y(r)}, \eqref{bc:FKPP:y(r)}.

\begin{example}\label{exam-z_0<z_1}\nopagebreak
\begingroup\rm
From a combination of \eqref{eq:FKPP:y(r)}, \eqref{bc:FKPP:y(r)},
\eqref{def:dz/dU}, and \eqref{lim:dz/dU}
one can deduce that $z_0 > -\infty$ occurs if the reaction term
$f(r)$ has the following asymptotic behavior as $r\to 0+$:
\begin{equation}
\label{e:f(0+)}
  \lim_{r\to 0+} \frac{f(r)}{r^{\alpha_0}} = {}- \gamma_0 < 0
  \quad\mbox{ where }\quad 0 < \alpha_0 < 1 \,.
\end{equation}
Analogously, $z_1 < +\infty$ occurs if
$f(r)$ has the following asymptotic behavior as $r\to 1-$:
\begin{equation}
\label{e:f(1-)}
  \lim_{r\to 1-} \frac{f(r)}{(1-r)^{\alpha_1}} = \gamma_1 > 0
  \quad\mbox{ where }\quad 0 < \alpha_1 < 1 \,.
\end{equation}
These claimes are proved in details in
{\sc P.\ Dr\'abek} and {\sc P.\ Tak\'a\v{c}} \cite{DrabTak-1},
Corrolary 6.3, {\rm Part~(i)}.
Notice that such a function $f$ cannot be differentiable or Lipschitzian
at the point $r=0$ or $r=1$, respectively.
\endgroup
\end{example}
\par\vskip 10pt

It is well\--known that in the classical case of
$f\colon \RR\to \RR$ being continuously differentiable, one has
$z_0 = -\infty$ and $z_1 = +\infty$, see
{\sc P.~C.\ Fife} and {\sc J.~B.\ Mc{L}eod}
\cite[Sect.~1]{Fife-McLeod} and
{\sc J.~D.\ Murray} \cite{Murray-I}, {\S}13.3, pp.\ 444--449.
At this point, we do not distinguish between the cases
$z_0$ or $z_1$ being finite or infinite.
Towards the end of this article, when investigating
the $\omega$-limit sets of solutions to
the initial value problem~\eqref{e_c:FKPP},
we will focus on the new case
$-\infty < z_0 < z_1 < \infty$ only,
which has not yet been treated in the literature.

\section{Solutions of the Initial Value Problem~\eqref{e_c:FKPP}}
\label{s:exist:FKPP}

In order to be able to establish the H\"older continuity of
all partial derivatives that appear in eq.~\eqref{e:FKPP} above,
we impose the following H\"older continuity hypothesis on
the reaction function $f$:
\begin{enumerate}
\renewcommand{\labelenumi}{{\bf (H\alph{enumi})}}
\item[{\bf (H2)}]
$f\colon \RR\to \RR$ is an $\alpha$-H\"older continuous function,
with an exponent $\alpha\in (0,1)$.
\end{enumerate}

We investigate the existence and uniqueness of a classical solution
$v(z,t)$ to the Cauchy (initial value) problem \eqref{e_c:FKPP}
obtained from problem \eqref{e:FKPP}.
We seek bounded ($L^{\infty}$-) classical solutions
``squeezed'' (i.e., bounded below and above, respectively)
between an ordered pair of a sub- and super\-solution to
the initial value problem~\eqref{e_c:FKPP}.
Such a solution takes values $v(z,t)\in [0,1]$
for all $z\in \RR$ and for all $t\in (0,\infty)$.
The initial condition is satisfied in the sense of
the weak$^{\ast}$\--limit in $L^{\infty}(\RR)$.
In order to obtain the desired regularity properties of
a {\em bounded classical solution\/}
to the initial value problem~\eqref{e_c:FKPP}, below we use
the Green function associated with the linear part of this problem.

\subsection{Mild and Classical Solutions}
\label{ss:Mild_Sol}

First, we denote by
\begin{equation}
\label{e:Green_heat}
\begin{aligned}
  G(x,y;t)\equiv G(|x-y|;t)
  \eqdef \frac{1}{ 2\sqrt{\pi t} }\,
  \exp\left( {}- \genfrac{}{}{}0{|x-y|^2}{4t} \right)
\\
  \quad\mbox{ for $x,y\in \RR$ and $t\in (0,\infty)$ }
\end{aligned}
\end{equation}
the Green function for the standard heat equation
(cf.\ \eqref{e:FKPP}).
The desired Green function for
the linear analogue of problem~\eqref{e_c:FKPP},
\begin{equation}
\label{e_c:lin-FKPP}
\left\{
\begin{aligned}
    \frac{\partial v}{\partial t}
  - \frac{\partial^2 v}{\partial z^2}
  - c\, \frac{\partial v}{\partial z}
& = g(z,t) \,,
    \quad (z,t)\in \RR\times \RR_+ \,;
\\
  v(z,0)&= v_0(z) \,,
\end{aligned}
\right.
\end{equation}
where we have replaced the nonlinearity $f(v(z,t))$
by a given function $g\in L^{\infty}(\RR\times \RR_+)$,
is obtained by shifting the space variable $z\mapsto x = z+ct$
in eq.~\eqref{e:Green_heat},
\begin{equation}
\label{e:Green_c}
\begin{aligned}
& G^{(c)}(z,y;t)\equiv G^{(c)}(z-y;t) \eqdef G(|z-y+ct|;t)
\\
& = \frac{1}{ 2\sqrt{\pi t} }\,
  \exp\left( {}- \genfrac{}{}{}0{|z-y+ct|^2}{4t} \right)
  \quad\mbox{ for $z,y\in \RR$ and $t\in (0,\infty)$. }
\end{aligned}
\end{equation}
Second, we require that any {\rm bounded classical solution\/}
$v\colon \RR\times \RR_+\to [0,1]\subset \RR$
to problem~\eqref{e_c:FKPP}
be also a {\em mild $L^{\infty}$\--solution\/}, i.e.,
$v$ must be essentially bounded and obey the following integral equation:
\begin{align}
\label{int_c:FKPP}
    v(z,t)
  = [ \mathcal{G}_1(t) v_0 ](z)
  + [ \mathcal{G}_2 (f\circ v) ](z,t)
  \quad\mbox{ for }\, (z,t)\in \RR\times (0,\infty) \,,
\end{align}
where $\mathcal{G}_1(t)$ and $\mathcal{G}_2$
are integral operators defined by
\begin{align}
\label{def:G_1}
    [ \mathcal{G}_1(t) g_0 ](z)
& \eqdef
    \int_{-\infty}^{\infty} G^{(c)}(z,y;t)\, g_0(y) \,\mathrm{d}y
  \quad\mbox{ and }\quad
\\
\label{def:G_2}
    [ \mathcal{G}_2 g ](z,t)
& \eqdef
    \int_0^t \int_{-\infty}^{\infty} G^{(c)}(z,y;t-s)\, g(y,s)
    \,\mathrm{d}y \,\mathrm{d}s
\end{align}
for $(z,t)\in \RR\times (0,\infty)$ and for all functions
$g_0\in L^{\infty}(\RR)$ and $g\in L^{\infty}(\RR\times \RR_+)$.

It is easy to see that the first operator,
$\mathcal{G}_1(t)$, has the following boundedness property:
Given any nonnegative integers $k,m\in \ZZ_+$,
there exists a constant $M_{k,m}\in \RR_+$ such that
\begin{align}
\label{est:G_1}
  \left\vert
  \frac{\partial^{k+m}}{\partial t^k\, \partial z^m}
    [ \mathcal{G}_1(t) g_0 ](z)
  \right\vert
  \leq M_{k,m}\, t^{-k-(m/2)}\cdot \| g_0\|_{ L^{\infty}(\RR) }
\\
\nonumber
  \quad\mbox{ holds for all }\, (z,t)\in \RR\times (0,\infty) \,.
\end{align}
This estimate follows directly from
{\sc O.~A.\ Ladyzhenskaya}, {\sc N.~N.\ Ural'tseva}, and
{\sc V.~A.\ Solonnikov}
\cite[Sect.~IV, {\S}1]{Ladyz-parab}, ineq.\ (2.5) on p.~274;
the partial derivative on the left\--hand side is taken
pointwise in the classical sense.
The second operator, $\mathcal{G}_2$, is a bit more complicated:
Given any positive number $T\in (0,\infty)$,
there exists a constant $M^{(T)}\in \RR_+$ such that
\begin{align}
\label{est_x:G_2}
  \left\vert
  \frac{\partial}{\partial z} [ \mathcal{G}_2 g ](z,t)
  \right\vert
  \leq M^{(T)}\, \| g\|_{ L^{\infty}(\RR\times (0,T)) }
    \quad\mbox{ and }\quad
\\
\label{est_t:G_2}
  \frac{ \left|
    [ \mathcal{G}_2 g ](z,t') - [ \mathcal{G}_2 g ](z,t)
         \right| }{ |t'- t|^{1/2} }
  \leq M^{(T)}\, \| g\|_{ L^{\infty}(\RR\times (0,T)) }
\\
\nonumber
  \quad\mbox{ hold for all }\,
  (z,t), (z,t')\in \RR\times [0,T] \,,\ t\neq t' \,.
\end{align}
Also these estimates follow directly from
\cite[Sect.~IV, {\S}1]{Ladyz-parab}, p.~263.

Finally, we quote the following regularity result for
$\mathcal{G}_2 g$ proved in
\cite[Sect.~IV, {\S}2]{Ladyz-parab}, ineq.\ (2.1) on p.~273:
Let $\ell\in (0,1)$ be arbitrary.
If $g\in C^{\ell, \ell/2}(\RR\times [0,T])$ then all partial derivatives
\begin{equation*}
  \frac{\partial}{\partial z} [ \mathcal{G}_2 g ] \,,\quad
  \frac{\partial^2}{\partial z^2} [ \mathcal{G}_2 g ] \,
    \quad\mbox{ and }\quad
  \frac{\partial}{\partial t} [ \mathcal{G}_2 g ]
\end{equation*}
belong to the H\"older space
$C^{\ell, \ell/2}(\RR\times [0,T])$ defined below.
Furthermore, there is a constant
$M_{\ell}^{(T)}\in \RR_+$, independent from
$g\in C^{\ell, \ell/2}(\RR\times [0,T])$, such that
\begin{align}
\nonumber
& \max\left\{
  \left\Vert
  \frac{\partial}{\partial z} [ \mathcal{G}_2 g ]
  \right\Vert_{ C^{\ell, \ell/2}(\RR\times [0,T]) } \,,\;
  \left\Vert
  \frac{\partial^2}{\partial z^2} [ \mathcal{G}_2 g ]
  \right\Vert_{ C^{\ell, \ell/2}(\RR\times [0,T]) } \,,\;
  \left\Vert
  \frac{\partial}{\partial t} [ \mathcal{G}_2 g ]
  \right\Vert_{ C^{\ell, \ell/2}(\RR\times [0,T]) }
    \right\}
\\
\label{Hoelder:G_2}
& \leq M_{\ell}^{(T)}\,
       \| g\|_{ C^{\ell, \ell/2}(\RR\times [0,T]) } \,.
\end{align}
This H\"older space is defined as follows.
Given any number $\ell\in (0,1)$, we denote by
$C^{\ell, \ell/2}(\RR\times [0,T])$
the Banach space of all bounded functions
$g\colon \RR\times [0,T]\to \RR$ such that
\begin{equation}
\label{e:Hoelder}
\begin{aligned}
  \left[ g\right]^{(\ell, \ell/2)}
& \eqdef
    \sup_{ 0 < |z'- z|\leq 1 ,\ t\in [0,T] }
  \frac{ |g(z',t) - g(z,t)| }{ |z'- z|^{\ell} }
\\
& {}
  + \sup_{ z\in \RR ,\ 0\leq t < t'\leq T }
  \frac{ |g(z,t') - g(z,t)| }{ |t'- t|^{\ell/2} }
  < \infty \,.
\end{aligned}
\end{equation}
The norm in $C^{\ell, \ell/2}(\RR\times [0,T])$ is defined by
\begin{equation}
\label{norm:Hoelder}
  \| g\|_{ C^{\ell, \ell/2}(\RR\times [0,T]) }
  \eqdef \left[ g\right]^{(\ell, \ell/2)}
  + \| g\|_{ L^{\infty}(\RR\times [0,T]) } \,.
\end{equation}

These regularity results for the integral operators
$\mathcal{G}_1(t)$ and $\mathcal{G}_2$ yield
the following differentiability properties for
any {\rm mild $L^{\infty}$\--solution\/}
$v$ of the integral equation \eqref{int_c:FKPP}:
Let $0 < t_0 < \tau < T < \infty$ and let us begin with
an arbitrary function
$g\in L^{\infty}(\RR\times (0,T))$ in place of $f\circ v$.
Then the function $v(z,t)$ on the left\--hand side of
eq.~\eqref{int_c:FKPP} satisfies
\begin{align}
\label{est_x:v_2}
& \left\vert
  \frac{\partial v}{\partial z} (z,t)
  \right\vert
  \leq C_{t_0}^{(T)}
    \left( \| v_0\|_{ L^{\infty}(\RR) }
  + \| g\|_{ L^{\infty}(\RR\times (0,T)) }
    \right)
    \quad\mbox{ and }\quad
\\
\label{est_t:v_2}
& \frac{ \left| v(z,t') - v(z,t) \right| }{ |t'- t|^{1/2} }
  \leq C_{t_0}^{(T)}
    \left( \| v_0\|_{ L^{\infty}(\RR) }
  + \| g\|_{ L^{\infty}(\RR\times (0,T)) }
    \right)
\\
\nonumber
& \quad\mbox{ for all }\,
  (z,t), (z,t')\in \RR\times [t_0,T] \,,\ t\neq t' \,,
\end{align}
where $C_{t_0}^{(T)}\in \RR_+$ is a constant independent from
$v_0\in L^{\infty}(\RR)$ and $g\in L^{\infty}(\RR\times (0,T))$.
These estimates follow from an easy combination of the inequalities in
\eqref{est:G_1}, \eqref{est_x:G_2}, and \eqref{est_t:G_2},
with the interval $[0,T]$ replaced by $[t_0,T]$.

The estimates below take advantage of the fact that
the linear initial value problem \eqref{e_c:lin-FKPP},
with any given inhomogeneity $g\in L^{\infty}(\RR\times \RR_+)$
in place of $f\circ v$ on the right\--hand side,
possesses a unique mild $L^{\infty}$\--solution
$v\colon \RR\times \RR_+\to \RR$.
In particular, this solution,
\begin{equation*}
    v(z,t)
  = [ \mathcal{G}_1(t) v_0 ](z)
  + [ \mathcal{G}_2 g ](z,t) \,,
  \quad\mbox{ for }\, (z,t)\in \RR\times (0,\infty) \,,
\end{equation*}
satisfies \eqref{est_x:v_2} and \eqref{est_t:v_2}.
Hence, recalling our {\rm Hypothesis {\bf (H2)}}
(i.e., $f\colon \RR\to \RR$ is $\alpha$-H\"older\--continuous)
and setting $\tilde{g} = f\circ v$,
\eqref{est_x:v_2} and \eqref{est_t:v_2} imply
$\tilde{g}\in C^{\alpha, \alpha/2}(\RR\times [t_0,T])$.
We denote by $\tilde{v}$ the unique mild $L^{\infty}$\--solution of
problem \eqref{e_c:lin-FKPP} in the domain $\RR\times (t_0,T)$,
\begin{equation*}
\left\{
\begin{aligned}
    \frac{\partial \tilde{v}}{\partial t}
  - \frac{\partial^2 \tilde{v}}{\partial z^2}
  - c\, \frac{\partial \tilde{v}}{\partial z}
& = \tilde{g} \,,
    \quad (z,t)\in \RR\times (t_0,T) \,;
\\
  \tilde{v}(z,t_0)&= v(z,t_0) \,.
\end{aligned}
\right.
\end{equation*}
Recall that the initial value
$v(\,\cdot\,,t_0)\in L^{\infty}(\RR)$ satisfies \eqref{est_x:v_2}.
Here, we have replaced the initial time $t=0$ by $t = t_0\in (0,\tau)$.
We combine the regularity estimates in
\eqref{est:G_1} and \eqref{Hoelder:G_2}, thus arriving at
\begin{align}
\nonumber
& \max\left\{
  \left\Vert \frac{\partial \tilde{v}}{\partial z}
  \right\Vert_{ C^{\alpha, \alpha/2}(\RR\times [\tau,T]) } \,,\;
  \left\Vert \frac{\partial^2 \tilde{v}}{\partial z^2}
  \right\Vert_{ C^{\alpha, \alpha/2}(\RR\times [\tau,T]) } \,,\;
  \left\Vert \frac{\partial \tilde{v}}{\partial t}
  \right\Vert_{ C^{\alpha, \alpha/2}(\RR\times [\tau,T]) }
    \right\}
\\
\label{Hoelder:^v_2}
& \leq C_{\alpha}^{(t_0,\tau,T)}
    \left( \| v(\,\cdot\,,t_0)\|_{ L^{\infty}(\RR) }
  + \|\tilde{g}\|_{ C^{\alpha, \alpha/2}(\RR\times [t_0,T]) }
    \right)
\\
\nonumber
& \leq C_{\alpha}^{(\tau,T)}
    \left( \| v\|_{ L^{\infty}(\RR\times (0,T)) }
  + \| g\|_{ L^{\infty}(\RR\times (0,T)) }
    \right) \,,
\end{align}
where
\begin{math}
  C_{\alpha}^{(t_0,\tau,T)}\in \RR_+
\end{math}
is a constant independent from
$v(\,\cdot\,,t_0)\in L^{\infty}(\RR)$ and
$\tilde{g}\in C^{\alpha, \alpha/2}(\RR\times [t_0,T])$, and
\begin{math}
  C_{\alpha}^{(\tau,T)}\in \RR_+
\end{math}
is another constant independent from
$v\in L^{\infty}(\RR\times (0,T))$ and
$g\in L^{\infty}(\RR\times (0,T))$.
To derive the second inequality in \eqref{Hoelder:^v_2},
we have made use of \eqref{est_x:v_2} and \eqref{est_t:v_2}.

We conclude from \eqref{Hoelder:^v_2} that
any mild $L^{\infty}$\--solution $v$ of problem \eqref{e_c:FKPP}
satisfies
\begin{equation*}
  \frac{\partial v}{\partial z} \,,\
  \frac{\partial^2 v}{\partial z^2} \,,\
  \frac{\partial v}{\partial t}
  \in C^{\alpha, \alpha/2}(\RR\times [\tau,T])
    \quad\mbox{ whenever }\, 0 < \tau < T < \infty \,,
\end{equation*}
and
\begin{align}
\label{Hoelder:v_2}
& \max\left\{
  \left\Vert \frac{\partial v}{\partial z}
  \right\Vert_{ C^{\alpha, \alpha/2}(\RR\times [\tau,T]) } \,,\;
  \left\Vert \frac{\partial^2 v}{\partial z^2}
  \right\Vert_{ C^{\alpha, \alpha/2}(\RR\times [\tau,T]) } \,,\;
  \left\Vert \frac{\partial v}{\partial t}
  \right\Vert_{ C^{\alpha, \alpha/2}(\RR\times [\tau,T]) }
    \right\}
\\
\nonumber
& \leq C_{\alpha}^{(\tau,T)}
    \left( \| v\|_{ L^{\infty}(\RR\times (0,T)) }
  + \| f\circ v\|_{ L^{\infty}(\RR\times (0,T)) }
    \right) \,.
\end{align}
Consequently, $v$ is also
a {\em bounded classical solution\/}
to problem~\eqref{e_c:FKPP};
the initial condition is satisfied in the sense of
the weak$^{\ast}$\--limit
$v(\,\cdot\,,t) \wstarconverge v_0$
in $L^{\infty}(\RR)$ as $t\to 0+$.

\subsection{Weak Sub- and Supersolutions and Weak Solutions}
\label{ss:sub-sup}

Let $v_0\in L^{\infty}(\RR)$.
A function
$\underline{v}\in L^{\infty}(\RR\times \RR_+)$
will be called a {\em weak $L^{\infty}$\--sub\-solution\/}
of the initial value problem \eqref{e_c:FKPP}
if it satisfies the following three conditions:
\begin{itemize}
\item[{\rm (i)}]
For every nonnegative test function $\phi\in W_0^{1,1}(\RR)$,
we have
\begin{align*}
  \limsup_{t\to 0+}
    \int_{-\infty}^{\infty}
    [ \underline{v}(z,t) - v_0(z) ]\, \phi(z) \,\mathrm{d}z
  \leq 0 \,.
\end{align*}
\item[{\rm (ii)}]
$\underline{v}$ is Lipschitz\--continuous in every set
$\RR\times [\tau,T]$ whenever $0 < \tau < T < \infty$, i.e.,
\begin{math}
  {\partial \underline{v}} / {\partial z} \,,\
  {\partial \underline{v}} / {\partial t}
  \in L^{\infty}(\RR\times [\tau,T]) \,.
\end{math}
\item[{\rm (iii)}]
For every nonnegative test function $\phi\in W_0^{1,1}(\RR)$,
the following inequality holds for a.e.\ $t\in (0,\infty)$,
\begin{align}
\label{sub:FKPP}
& \frac{\mathrm{d}}{\mathrm{d}t}
    \int_{-\infty}^{\infty} \underline{v}(z,t)\, \phi(z) \,\mathrm{d}z
  + \int_{-\infty}^{\infty}
    \frac{ \partial\underline{v} }{\partial z}(z,t)\cdot
    \frac{\partial\phi}{\partial z}(z) \,\mathrm{d}z
\\
\nonumber
& {}
  - c \int_{-\infty}^{\infty}
    \frac{ \partial\underline{v} }{\partial z}(z,t)\cdot \phi(z)
    \,\mathrm{d}z
  \leq \int_{-\infty}^{\infty}
    f(\underline{v}(z,t))\, \phi(z) \,\mathrm{d}z \,.
\end{align}
\end{itemize}

A {\em weak $L^{\infty}$\--super\-solution\/} 
$\overline{v}\in L^{\infty}(\RR\times \RR_+)$
is defined analogously.
We say that a function
$v\in L^{\infty}(\RR\times \RR_+)$
is a {\em weak $L^{\infty}$\--solution\/}
of the initial value problem \eqref{e_c:FKPP}
if and only if
$v$ is a weak $L^{\infty}$\--sub- and -super\-solution
to \eqref{e_c:FKPP}.
Of course, in Condition {\rm (i)},
the nonnegative test functions $\phi\in W_0^{1,1}(\RR)$
may be replaced by nonnegative functions $\phi\in L^1(\RR)$.
Condition {\rm (ii)} combined with equality in
\eqref{sub:FKPP} implies
\begin{math}
  {\partial v} / {\partial z} \,,\
  {\partial v} / {\partial t} \,,\
  {\partial^2 v} / {\partial z^2}
  \in L^{\infty}(\RR\times [\tau,T]) \,.
\end{math}

We remark that any weak $L^{\infty}$\--solution of
problem \eqref{e_c:FKPP}
is also a mild $L^{\infty}$\--solution and vice versa; cf.\
{\sc J.~M.\ Ball} \cite[Theorem, p.~371]{Ball} or
{\sc A.\ Pazy} \cite[Theorem, p.~259]{Pazy}.
Thanks to the regularity properties in \eqref{Hoelder:v_2},
a mild $L^{\infty}$\--solution is also a bounded classical solution.
It is shown in
{\sc A.\ Pazy} \cite[{\S}4.2, pp.\ 105--110]{Pazy}
that any bounded classical solution is also
a mild $L^{\infty}$\--solution.

All types of sub- and super\-solutions and solutions of
problem \eqref{e_c:FKPP}, defined above for all times $t\in \RR_+$,
can be defined analogously for time $t\in [\tau,T)$
from a bounded time interval $[\tau,T)\subset \RR_+$.

The following {\it one\--sided Lipschitz condition\/}
is a crucial hypothesis imposed on the function $f$
in our method for establishing the {\em weak comparison principle\/}
for weak $L^{\infty}$-sub- and -super\-solutions:
\begin{enumerate}
\renewcommand{\labelenumi}{{\bf (H\alph{enumi})}}
\item[{\bf (H3)}]
There is a number $L\in \RR_+$ such that
\begin{equation}
\label{e:f_Lip}
  f(s') - f(s)\leq L\, (s'- s)
    \quad\mbox{ for all }\, s,s'\in \RR \,,\ s < s' \,.
\end{equation}
\end{enumerate}
%

\begin{lemma}\label{lem-weak_comp}
{\rm (Weak comparison principle.)}$\,$
Let\/ $c\in \RR$ and let\/
$f\colon \RR\to \RR$ be continuous and satisfy\/
{\rm Hypothesis {\bf (H3)}}.
Assume that\/
$\underline{v}, \overline{v}\in L^{\infty}(\RR\times \RR_+)$
is a pair of weak $L^{\infty}$-sub- and -super\-solutions
to problem~\eqref{e_c:FKPP}, such that\/
\begin{math}
  \underline{v}(\,\cdot\,,0) = \underline{v}_0
  \leq \overline{v}_0 = \overline{v}(\,\cdot\,,0)
\end{math}
in $L^{\infty}(\RR)$.
Then we have also $\underline{v}\leq \overline{v}$ a.e.\ in
$\RR\times \RR_+$.
\end{lemma}

\par\vskip 10pt
\proof
We subtract the analogue of ineq.\ \eqref{sub:FKPP}
for a supersolution (with the reversed inequality)
from \eqref{sub:FKPP} for a subsolution,
\begin{align}
\nonumber
& \frac{\mathrm{d}}{\mathrm{d}t}
    \int_{-\infty}^{\infty}
    ( \underline{v}(z,t) - \overline{v}(z,t) )\, \phi(z) \,\mathrm{d}z
  + \int_{-\infty}^{\infty}
    \frac{\partial}{\partial z}
    ( \underline{v} - \overline{v} )\cdot
    \frac{\partial\phi}{\partial z}(z) \,\mathrm{d}z
\\
\label{e_c:sub_v-sup_v}
& {}
  - c \int_{-\infty}^{\infty}
    \frac{\partial}{\partial z}
    ( \underline{v} - \overline{v} )\cdot \phi(z) \,\mathrm{d}z
\\
\nonumber
& \leq \int_{-\infty}^{\infty}
    [ f( U(z) + \underline{v}(z) ) - f( U(z) + \overline{v}(z) ) ]\,
    \phi(z) \,\mathrm{d}z \,,
\end{align}
for every nonnegative test function $\phi\in W_0^{1,1}(\RR)$.
Recalling Conditions {\rm (i)} and {\rm (ii)}, we observe that
\begin{align*}
\nonumber
  \frac{\mathrm{d}}{\mathrm{d}t}
    \int_{-\infty}^{\infty}
    ( \underline{v}(z,t) - \overline{v}(z,t) )\, \phi(z) \,\mathrm{d}z
  = \int_{-\infty}^{\infty}
    \frac{\partial}{\partial t}
    ( \underline{v} - \overline{v} )\cdot \phi(z) \,\mathrm{d}z
\end{align*}
holds for a.e.\ $t\in (0,\infty)$.

Let
\begin{math}
    (\underline{v} - \overline{v})^{+}
  = \max\{ (\underline{v} - \overline{v}) ,\, 0\}
\end{math}
denote the nonnegative part of the difference
$\underline{v} - \overline{v}$.
By standard arguments, we have
$(\underline{v} - \overline{v})^{+}\in L^{\infty}(\RR\times \RR_+)$
and
\begin{equation*}
  \frac{\partial}{\partial t}
    (\underline{v} - \overline{v})^{+} \,,\
  \frac{\partial}{\partial z}
    (\underline{v} - \overline{v})^{+}
  \in L^{\infty}(\RR\times [\tau,T])
  \quad\mbox{ whenever }\, 0 < \tau < T < \infty \,.
\end{equation*}
For a.e.\ $t\in (0,\infty)$, in ineq.~\eqref{e_c:sub_v-sup_v} above
we may replace $\phi$ by the product
$(\underline{v} - \overline{v})^{+}(\,\cdot\,,t)\, \phi$,
thus obtaining
\begin{align*}
\nonumber
& \frac{1}{2}\cdot \frac{\mathrm{d}}{\mathrm{d}t}
    \int_{-\infty}^{\infty}
  [ (\underline{v} - \overline{v})^{+}(z,t) ]^2\cdot \phi(z) \,\mathrm{d}z
\\
& {}
  + \int_{-\infty}^{\infty}
    \frac{\partial}{\partial z}
    ( \underline{v} - \overline{v} )\cdot
    \frac{\partial}{\partial z}
  \left[
    (\underline{v} - \overline{v})^{+}(z,t)\, \phi(z)
  \right] \,\mathrm{d}z
\\
\nonumber
& {}
  - c \int_{-\infty}^{\infty}
    \frac{\partial}{\partial z}
    ( \underline{v} - \overline{v} )\cdot
    (\underline{v} - \overline{v})^{+}(z,t)\, \phi(z) \,\mathrm{d}z
\\
\nonumber
& \leq \int_{-\infty}^{\infty}
    [ f( U(z) + \underline{v}(z) ) - f( U(z) + \overline{v}(z) ) ]\,
  (\underline{v} - \overline{v})^{+}(z,t)\, \phi(z) \,\mathrm{d}z \,.
\end{align*}
We simplify the integrands and apply inequality \eqref{e:f_Lip}
from the one\--sided Lipschitz condition in {\rm Hypothesis {\bf (H3)}}
to the last integral, thus arriving at
\begin{align}
\nonumber
& \frac{1}{2}\cdot \frac{\mathrm{d}}{\mathrm{d}t}
    \int_{-\infty}^{\infty}
  [ (\underline{v} - \overline{v})^{+}(z,t) ]^2\cdot \phi(z) \,\mathrm{d}z
\\
\label{ineq:sub_v-sup_v}
& {}
  + \frac{1}{2}
    \int_{-\infty}^{\infty}
    \frac{\partial}{\partial z}
  [ (\underline{v} - \overline{v})^{+}(z,t) ]^2\cdot
    \frac{\partial\phi}{\partial z}(z) \,\mathrm{d}z
  + \int_{-\infty}^{\infty}
  \left[ \frac{\partial}{\partial z} (\underline{v} - \overline{v})^{+}
  \right]^2\cdot \phi(z) \,\mathrm{d}z
\\
\nonumber
& {}
  - \frac{c}{2} \int_{-\infty}^{\infty}
    \frac{\partial}{\partial z}
  [ (\underline{v} - \overline{v})^{+}(z,t) ]^2\cdot
    \phi(z) \,\mathrm{d}z
\\
\nonumber
& \leq L\int_{-\infty}^{\infty}
  [ (\underline{v} - \overline{v})^{+}(z,t) ]^2\cdot
    \phi(z) \,\mathrm{d}z \,.
\end{align}

We conclude that the nonnegative function
$W(z,t) = \ee^{2Lt}\, [ (\underline{v} - \overline{v})^{+}(z,t) ]^2$
satisfies the inequality
\begin{equation}
\label{e:sub_v-sup_v}
\left\{
\begin{aligned}
    \frac{\partial W}{\partial t}
  - \frac{\partial^2 W}{\partial z^2}
  - c\, \frac{\partial W}{\partial z}
& \leq 0 \,,
    \quad (z,t)\in \RR\times \RR_+ \,;
\\
  W(z,0)&= 0 \,,
\end{aligned}
\right.
\end{equation}
in the weak sense,
$W\in L^{\infty}(\RR\times (0,T))$ for every $T\in (0,\infty)$.
Consequently, we substitute $x = z+ct$ and apply
the weak maximum principle for the heat equation to the function
$W(z+ct,t)$, which yields
$W\leq 0$ a.e.\ in $\RR\times \RR_+$; hence,
$(\underline{v} - \overline{v})^{+} \equiv 0$
a.e.\ in $\RR\times \RR_+$.

The lemma is proved.
\qed
\par\vskip 10pt

Lemma~\ref{lem-weak_comp}
has the following straight\--forward consequence.

\begin{corollary}\label{cor-weak_comp}
{\rm (Uniqueness.)}$\,$
Let\/ $c\in \RR$ and let\/
$f\colon \RR\to \RR$ be continuous and satisfy\/
{\rm Hypothesis {\bf (H3)}}.
Then the initial value problem \eqref{e_c:FKPP}
has {\rm at most one\/} weak $L^{\infty}$\--solution.
In particular, the same uniqueness result holds also for
a bounded classical solution, as well.
Finally, if also {\rm Hypothesis {\bf (H2)}} is satisfied,
then this uniqueness result applies also to
any mild $L^{\infty}$\--solution.
\end{corollary}
\par\vskip 10pt

We give the proof of the existence of a weak $L^{\infty}$\--solution
to the initial value problem \eqref{e_c:FKPP}
in the next paragraph ({\S}\ref{ss:exist:FKPP}).

\subsection{Existence of a Solution to Problem \eqref{e_c:FKPP}}
\label{ss:exist:FKPP}

The following existence result complements our uniqueness result from
Corollary~\ref{cor-weak_comp}.

\begin{proposition}\label{prop-weak_sol}
Let\/ $c\in \RR$ and let\/
$f\colon \RR\to \RR$ satisfy\/
{\rm Hypotheses} {\bf (H2)} and {\bf (H3)}.
Assume that the initial data $v_0\colon \RR\to \RR$
are Lebesgue\--measurable and satisfy\/
$0\leq v_0\leq 1$ a.e.\ on~$\RR$.
Then the initial value problem \eqref{e_c:FKPP}
possesses {\rm a unique\/} mild $L^{\infty}$\--solution, say,
$v\colon \RR\times \RR_+\to \RR$.
In particular, the same existence and uniqueness result holds also for
weak $L^{\infty}$\--solutions and bounded classical solutions, as well.
Finally, $0\leq v(z,t)\leq 1$ holds for all\/
$(z,t)\in \RR\times (0,\infty)$.
\end{proposition}

\par\vskip 10pt
\proof
We apply the well\--known Tikhonov fixed point theorem
({\sc K.\ Deimling} \cite[Theorem 10.1, p.~90]{Deimling})
as follows.
Given any $T\in (0,\infty)$, we denote by
\begin{equation*}
  \mathcal{X} =
  C_{\mathrm{b}}(\RR\times [0,T]) \eqdef
  L^{\infty}(\RR\times (0,T)) \cap C(\RR\times [0,T])
\end{equation*}
the vector space of all bounded continuous functions
$v\colon \RR\times [0,T]\to \RR$
endowed with the locally convex topology of uniform convergence
on every compact set
$K = J\times [0,T]\subset \RR\times [0,T]$, where
$J = [a,b]\subset \RR$ is a compact interval.
Thus, $\mathcal{X}$ is a Fr\'echet space.

By our definition of a mild $L^{\infty}$\--solution
$v\colon \RR\times \RR_+\to [0,1]\subset \RR$
to problem~\eqref{e_c:FKPP},
$v$ must be essentially bounded and obey the integral equation
\eqref{int_c:FKPP}.
Accordingly, we split it as
\begin{equation}
\label{e:v=G_1+w}
  v(z,t) = [ \mathcal{G}_1(t) v_0 ](z) + w(z,t)
  \quad\mbox{ for }\, (z,t)\in \RR\times [0,T] \,,
\end{equation}
where $w\in \mathcal{X}$ is a fixed point of the self\--mapping
\begin{math}
  \mathcal{T}\colon \mathcal{X}\to \mathcal{X}
\end{math}
defined by
\begin{align}
\label{def:int_c:FKPP}
  (\mathcal{T} w)(z,t) \eqdef
  [ \mathcal{G}_2 (f\circ v) ](z,t)
  \quad\mbox{ for }\, (z,t)\in \RR\times [0,T] \,,
\\
\nonumber
  \mbox{ with $v$ defined in eq.~\eqref{e:v=G_1+w}, }\,
         v\in L^{\infty}(\RR\times (0,T)) \,.
\end{align}
We recall that $\mathcal{G}_1(t)$ and $\mathcal{G}_2$
are integral operators defined by eqs.\
\eqref{def:G_1} and \eqref{def:G_2}, respectively.
Next, let us consider the closed convex subset
\begin{equation*}
  \mathcal{C} =
  \{ w\in \mathcal{X}\colon |w(z,t)|\leq 1
     \,\mbox{ for all }\, (z,t)\in \RR\times [0,T] \}
\end{equation*}
of $\mathcal{X}$ and denote
\begin{equation*}
  M = \max_{-1\leq s\leq 2} |f(s)| \,;\quad 0 < M < \infty \,.
\end{equation*}
Recall that $0\leq v_0\leq 1$ on $\RR$ and the kernel of
the integral operators $\mathcal{G}_1(t)$ and $\mathcal{G}_2$ satisfies
$G^{(c)}(z,y;t) > 0$ and
$\int_{-\infty}^{\infty} G^{(c)}(z,y;t) \,\mathrm{d}y = 1$.
Consequently, for every $w\in \mathcal{C}$ we have
\begin{align*}
& 0 \leq [ \mathcal{G}_1(t) v_0 ](z)\leq 1 \,,\quad
  -1\leq v(z,t) = [ \mathcal{G}_1(t) v_0 ](z) + w(z,t)
    \leq 2 \,,
\\
& \quad\mbox{ and }\quad
  |f(v(z,t))|\leq M
  \quad\mbox{ for }\, (z,t)\in \RR\times [0,T] \,.
\end{align*}
Applying these inequalities to eq.~\eqref{def:int_c:FKPP},
we arrive at
\begin{align*}
  | (\mathcal{T} w)(z,t) | =
  \left| [ \mathcal{G}_2 (f\circ v) ](z,t) \right| \leq Mt
  \quad\mbox{ for }\, (z,t)\in \RR\times [0,T]
\end{align*}
and for every function $w\in \mathcal{C}$.
We take $T = 1/M\in (0,\infty)$ and observe that
$\mathcal{T}$ maps $\mathcal{C}$ into itself.

From the uniform continuity of $f$ on the compact interval
$[-1,2]$ and the properties of the kernel
$G^{(c)}(z,y;t-s)$ of the integral operator $\mathcal{G}_2$
we deduce that
$\mathcal{T}\colon \mathcal{C}\to \mathcal{C}$
is continuous.
Finally, we combine the regularity estimates
\eqref{est_x:G_2} and \eqref{est_t:G_2}
with Arzel\`a\--Ascoli's compacness criterion in $\mathcal{X}$
({\sc P.\ Dr\'abek} and {\sc J.\ Milota}
 \cite[Theorem 1.2.13, p.~32]{DrabMil})
to conclude that the image $\mathcal{T}(\mathcal{C})$
is relatively compact.
By Tikhonov's fixed point theorem,
$\mathcal{T}$ has a fixed point in $\mathcal{C}$, say,
$\hat{w}\in \mathcal{C}$.
The sum corresponding to eq.~\eqref{e:v=G_1+w},
\begin{align*}
    \hat{v}(z,t)
& = [ \mathcal{G}_1(t) v_0 ](z) + \hat{w}(z,t)
\\
& = [ \mathcal{G}_1(t) v_0 ](z)
  + [ \mathcal{G}_2 (f\circ \hat{v}) ](z,t)
  \quad\mbox{ for }\, (z,t)\in \RR\times [0,T] \,,
\end{align*}
is a mild $L^{\infty}$\--solution to problem~\eqref{e_c:FKPP}
on the bounded time interval $[0,T]$ in place of $\RR_+$
and, by regularity results in {\S}\ref{ss:Mild_Sol},
also a weak $L^{\infty}$\--solution and a bounded classical solution.
By an analogue of Lemma~\ref{lem-weak_comp} (weak comparison principle)
for the bounded time interval $[0,T]$, we have
$0\leq \hat{v}(z,t)\leq 1$ for all $(z,t)\in \RR\times [0,T]$.
As in the case of the weak comparison principle,
the uniqueness of $\hat{v}$ follows from an analogue of
Corollary~\ref{cor-weak_comp}.

Repeating this procedure in every time interval
$[t_0, t_0 + T]$ of length $T$, for each $t_0\in \RR_+$,
we can construct a mild $L^{\infty}$\--solution
to problem~\eqref{e_c:FKPP} for all times $t\in \RR_+$.
This is also a weak $L^{\infty}$\--solution and
a bounded classical solution.
Finally, the uniqueness follows from Corollary~\ref{cor-weak_comp}.
\qed
\par\vskip 10pt

\subsection{Sub- and Supersolutions Resulting from Travelling Waves}
\label{ss:Trav-sub/super}

Of course, the constant functions
$\underline{v}\equiv 0$ and $\overline{v}\equiv 1$
form a trivial ordered pair of
weak $L^{\infty}$-sub- and -super\-solutions
to problem \eqref{e_c:FKPP}.
Now we are ready to modify the TW\--solutions in order to construct
suitable nontrivial
weak $L^{\infty}$-sub- and -super\-solutions
to problem \eqref{e_c:FKPP},
as defined in {\S}\ref{ss:sub-sup}.

In addition to $f\colon \RR\to \RR$ being continuous,
we assume that $f$ satisfies also
the following two {\it ``secant'' conditions\/}
on the interval $[0,1]$:
\begin{enumerate}
\renewcommand{\labelenumi}{{\bf (H\alph{enumi})}}
\item[{\bf (H4)}]
There exists a number $\eta_0\in \RR$,
\begin{equation}
\label{e:delta<1/3}
  0 < \eta_0 < \frac{1}{3}\cdot \min\{ s_0 ,\, 1 - s_0\} \,,
\end{equation}
with the following property:
Given any $\eta\in (0, \eta_0]$,
there are constants $\delta\in (0,\eta)$ and
$\underline{\mu} > 0$, $\overline{\mu} > 0$, depending on $\eta$,
such that the following two conditions hold:
\begin{alignat}{2}
\label{ineq:f(0)}
  \inf_{ 0\leq s\leq \delta } [ f(s) - f(s+q) ]
& \geq \underline{\mu}\, q
&&  \quad\mbox{ for all }\, q\in ( 0, s_0 - 2\eta ] \,,
\\
\label{ineq:f(1)}
  \inf_{ 1 - \delta\leq s\leq 1 } [ f(s-q) - f(s) ]
& \geq \overline{\mu}\, q
&&  \quad\mbox{ for all }\, q\in (0, 1 - s_0 - 2\eta ] \,.
\end{alignat}
\end{enumerate}
In particular, a continuously differentiable function $f$
(the classical case \cite{Fife-McLeod}, \cite[{\S}13.3]{Murray-I})
satisfies {\bf (H4)} whenever
$f'(0) < 0$, $f'(s_0) > 0$, and $f'(1) < 0$.
However, the secant conditions are satisfied
under somewhat different hypotheses on the differentiability of $f$:

\begin{example}\label{exam-secant}\nopagebreak
\begingroup\rm
Our {\rm Hypothesis} {\bf (H4)} above is satisfied if
$f\colon \RR\to \RR$ is a continuous function such that
$f(0) = f(s_0) = f(1) = 0$ for some $0 < s_0 < 1$, together with
$f(s) < 0$ for every $s\in (0,s_0)$,
$f(s) > 0$ for every $s\in (s_0,1)$,
and $f$ is differentiable in
$(0, \eta^{*})\cup \{ s_0\}\cup (1 - \eta^{*}, 1)$
for some $\eta^{*}\in \RR$,
\begin{equation*}
  0 < \eta^{*} < \frac{1}{3}\cdot \min\{ s_0 ,\, 1 - s_0\} \,,
\end{equation*}
with the derivatives satisfying $f'(s_0) > 0$ and
\begin{equation}
\label{e:f'(u)}
     \lim_{s\to 0+} f'(s)
   = \lim_{s\to 1-} f'(s) = -\infty \,.
\end{equation}
The last two conditions guarantee inequalities
\eqref{ineq:f(0)} and \eqref{ineq:f(1)}, respectively,
for $q > 0$ small enough, whereas condition $f'(s_0) > 0$
guarantees them for
$q\in ( 0, s_0 - 2\eta ]$ near $s_0 - 2\eta$ and for
$q\in ( 0, 1 - s_0 - 2\eta ]$ near $1 - s_0 - 2\eta$, respectively,
if $\eta > 0$ is small enough, say,
\begin{math}
  0 < \eta\leq \eta_0\leq \eta^{*} \,.
\end{math}
If $\eta_0\leq q\leq s_0 - 2\eta_0$ then ineq.~\eqref{ineq:f(0)}
follows from a combination of conditions
$\lim_{s\to 0+} f'(s) = -\infty$ and $f'(s_0) > 0$
with $f(s) < 0$ for every $s\in (0,s_0)$.
In this case $s+q$ satisfies
$\eta_0\leq s+q\leq s_0 - 2\eta_0 + \delta < s_0 - \eta_0$,
thanks to
$0\leq s\leq \delta < \eta\leq \eta_0$.
On the other hand,
if $\eta_0\leq q\leq 1 - s_0 - 2\eta_0$ then ineq.~\eqref{ineq:f(1)}
follows from a combination of conditions
$\lim_{s\to 1-} f'(s) = -\infty$ and $f'(s_0) > 0$
with $f(s) > 0$ for every $s\in (s_0,1)$.
In this case $s-q$ satisfies
$s_0 + \eta_0 < s-q\leq 1 - \eta_0$, by
\begin{align*}
    s_0 + \eta_0
& {}
  = (1 - \eta_0) - (1 - s_0 - 2\eta_0)
  \leq (1 - \eta) - (1 - s_0 - 2\eta_0)
\\
& {}
  < (1 - \delta) - (1 - s_0 - 2\eta_0)
  \leq s-q\leq 1 - \eta_0 \,,
\end{align*}
thanks to
$1 - \eta_0\leq 1 - \eta < 1 - \delta\leq s\leq 1$.
\endgroup
\end{example}
\par\vskip 10pt

Finally, we assume:
\begin{enumerate}
\renewcommand{\labelenumi}{{\bf (H\alph{enumi})}}
\item[{\bf (H5)}]
The initial condition $v_0\in L^{\infty}(\RR)$
in the Cauchy problem \eqref{e_c:FKPP}
is defined at every point $z\in \RR$ and satisfies
$0\leq v_0(z)\leq 1$ together with
\begin{equation}
\label{e:v_0:sub-sup}
  v_0(-\infty)\eqdef
    \limsup_{z\to -\infty} v_0(z) < s_0
  < \liminf_{z\to +\infty} v_0(z)
  \eqdef v_0(+\infty) \,.
\end{equation}
\end{enumerate}

Accordingly, the constant $\eta\in (0,\eta_0]$
in {\rm Hypothesis} {\bf (H4)} is chosen to be small enough,
such that also
\begin{equation}
\label{eq:v_0:sub-sup}
  v_0(-\infty)\eqdef
    \limsup_{z\to -\infty} v_0(z) < s_0 - 2\eta < s_0 + 2\eta
  < \liminf_{z\to +\infty} v_0(z)
  \eqdef v_0(+\infty) \,.
\end{equation}
As in Theorem~\ref{thm-conv_TW}, we mean the monotone limits
\begin{equation*}
  v_0(-\infty) = \lim_{n\to \infty}\;
    \esssup_{-\infty < z\leq -n} v_0(z) \quad\mbox{ and }\quad
  v_0(+\infty) = \lim_{n\to \infty}\;
    \essinf_{n\leq z < +\infty}  v_0(z) \,.
\end{equation*}

Following the proof of Lemma 4.1 in
{\sc P.~C.\ Fife} and {\sc J.~B.\ Mc{L}eod}
\cite[pp.\ 347--348]{Fife-McLeod},
we construct an ordered pair of weak $L^{\infty}$-sub- and -supersolutions
to eq.~\eqref{e_c:FKPP},
$v_1(z,t)\leq v_2(z,t)$ and
$v_1(z,0)\leq v_0(z)\leq v_2(z,0)$,
having the special forms rendered by the TW\--solution~$U$.
Recall that
$0 < U < 1$ and $U'> 0$ in $(z_0,z_1)\subset \RR$ together with
$\lim_{z\to z_0+} U(z) = 0$ and
$\lim_{z\to z_1-} U(z) = 1$.

\begin{proposition}\label{prop-sub-sup}
Assume that\/ $f$ satisfies {\rm Hypotheses} {\bf (H1)} and {\bf (H4)},
and\/ $v_0$ satisfies {\bf (H5)}.
Let\/
$U\colon \RR\to \RR$ denote the stationary solution of
eq.~\eqref{e_c:FKPP} described in
{\rm Proposition~\ref{prop-uniq_TW}, Part~{\rm (c)}}.
Then there exist some constants\/
$\mu, \nu, q_{0,i}\in (0,\infty)$ and\/
$\xi_{\infty,i}\in \RR$ such that the functions
\begin{align}
\label{e:subsol}
  v_1(z,t)
& = \max\{ U(z - \xi_1(t)) - q_1(t) ,\, 0\} \quad\mbox{ and }
\\
\label{e:supsol}
  v_2(z,t)
& = \min\{ U(z + \xi_2(t)) + q_2(t) ,\, 1\}
    \quad\mbox{ for }\, (z,t)\in \RR\times \RR_+ \,,
\end{align}
where
\begin{equation}
\label{e:q_i,xi_i}
  q_i(t) = q_{0,i}\, \ee^{-\mu t} \quad\mbox{ and }\quad
  \xi_i(t) = \xi_{\infty,i} - \nu\, q_i(t)
    \quad\mbox{ for }\, t\geq 0 \,;\quad i=1,2 \,,
\end{equation}
are weak $L^{\infty}$-sub- and -super\-solutions of\/
the Cauchy problem \eqref{e_c:FKPP}, respectively.
Finally, we have
$0\leq v_1(z,t)\leq v_2(z,t)\leq 1$ for all\/
$(z,t)\in \RR\times \RR_+$, together with
$v_1(z,0)\leq v_0(z)\leq v_2(z,0)$ for all\/ $z\in \RR$ at $t=0$.
\end{proposition}

\par\vskip 10pt
\proof
To begin with, we look for
some suitable continuously differentiable functions
$q_i\colon \RR_+\to (0,\infty)$ and
$\xi_i\colon \RR_+\to \RR$ with the limits
$\xi_i(t)\to \xi_{\infty,i}\in \RR$ and
$q_i(t)\to 0$ as $t\to \infty$; $i=1,2$.

As both cases of weak $L^{\infty}$-sub- and -super\-solutions,
$v_1$ and $v_2$, respectively, are similar,
we treat only the former one, the subsolution to eq.~\eqref{e_c:FKPP},
\begin{equation}
\label{e:v_1}
  v_1(z,t) = v(z,t) = \max\{ U(z - \xi(t)) - q(t) ,\, 0\} \,.
\end{equation}
We leave out the index $i=1$
in $\xi_1(t) = \xi(t)$ and $q_1(t) = q(t)$.
If no confusion is likely, we also suppress the argument in
$U = U( z - \xi(t) )$.

Recalling {\rm Hypothesis} {\bf (H5)} on $v_0$,
inequalities \eqref{e:v_0:sub-sup},
we first choose the constant $\eta\in (0,\eta_0]$
small enough, such that also
inequalities \eqref{eq:v_0:sub-sup} are valid; hence,
$s_0 + 2\eta < v_0(+\infty)$.
Of course, also
\eqref{e:delta<1/3} holds with $\eta$ in place of $\eta_0$.
Then we choose $q_0\in \RR$ such that
\begin{equation*}
  s_0 + 2\eta < 1 - q_0 < \min\{ v_0(+\infty) ,\, 1 - \eta\} \,,
\end{equation*}
that is,
\begin{equation}
\label{e:q_0}
    \max\{ 1 - v_0(+\infty) ,\, \eta\}
  < q_0 < 1 - s_0 - 2\eta \,.
\end{equation}
Hence, there is some $z'> 0$ large enough, such that
$U(z) - q_0\leq 1 - q_0\leq v_0(z + z')$ for all $z\geq 0$,
which yields
$U(z - z') - q_0\leq v_0(z)$ for all $z\geq z'$.
Since $\lim_{z\to -\infty} U(z) = 0$,
there is another $z''> 0$ large enough, such that
$U(z - z'') - q_0\leq 0\leq v_0(z)$ for all $z\leq z'$.
Setting
$z^{\ast} = \max\{ z',\, z''\} > 0$
and using the fact that $U\colon \RR\to [0,1]$ is monotone increasing,
we arrive at
$U(z - z^{\ast}) - q_0\leq v_0(z)$ for all $z\in \RR$.
This inequality shows that the initial condition for the subsolution
$v(z,t)$ is satisfied provided $\xi(0)\geq z^{\ast}$.

Now let us take the corresponding numbers $\delta\in (0,\eta)$ and
$\underline{\mu} > 0$, $\overline{\mu} > 0$
in {\rm Hypothesis} {\bf (H4)},
all of them depending on $\eta$ fixed above.
We need to distinguish among the following three cases,
$0\leq U\leq \delta$, $\delta\leq U\leq 1-\delta$, and
$1-\delta\leq U\leq 1$.
We write the underlying domain $\RR\times \RR_+$ as the union
\begin{math}
  \RR\times \RR_+ =
  \Omega_{\delta}^{(1)}\cup \Omega_{\delta}^{(2)}\cup
  \Omega_{\delta}^{(3)}
\end{math}
of the corresponding subsets
\begin{align}
\label{e:Omega_1}
  \Omega_{\delta}^{(1)} &\eqdef
  \{ (z,t)\in \RR\times \RR_+\colon
     0\leq U( z - \xi(t) ) \leq \delta\}
  \subset \RR\times \RR_+ \,,
\\
\label{e:Omega_2}
  \Omega_{\delta}^{(2)} &\eqdef
  \{ (z,t)\in \RR\times \RR_+\colon
     \delta\leq U( z - \xi(t) ) \leq 1 - \delta\}
  \subset \RR\times \RR_+ \,, \quad\mbox{ and }
\\
\label{e:Omega_3}
  \Omega_{\delta}^{(3)} &\eqdef
  \{ (z,t)\in \RR\times \RR_+\colon
     1 - \delta\leq U( z - \xi(t) ) \leq 1\}
  \subset \RR\times \RR_+ \,.
\end{align}

We begin with the \underline{third case}, $1-\delta\leq U\leq 1$:
Thanks to our choice of the function
$q\colon \RR_+\to \RR$, the inequalities
$0 < q(t)\leq q_0$ hold for all $t\geq 0$.
Hence, by ineq.~\eqref{e:q_0}, we have also
\begin{equation*}
    U(z - \xi(t)) - q(t)\geq U - q_0
  \geq (1 - \delta) - q_0 > (1 - \delta) - (1 - s_0 - 2\delta)
  = s_0 + \delta > 0
\end{equation*}
which shows that the subsolution $v_1 = v$ in eq.~\eqref{e:v_1}
is given by
\begin{equation}
\label{e:v_1:3}
  v(z,t) = U(z - \xi(t)) - q(t) \quad\mbox{ for all }\,
    (z,t)\in \Omega_{\delta}^{(3)} \,.
\end{equation}
It remains to verify the inequality
\begin{equation}
\label{ineq:sub}
    \mathcal{N}(v)\eqdef
    \frac{\partial v}{\partial t}
{}- \frac{\partial^2 v}{\partial z^2}
  - c\, \frac{\partial v}{\partial z} - f(v)
  \leq 0 \,,
    \quad (z,t)\in \Omega_{\delta}^{(3)} \,,
\end{equation}
which has to hold for a subsolution.
Using \eqref{e:dU/dz} and \eqref{e:v_1:3} we calculate
\begin{align}
\label{e:dv/dt}
&
    \frac{\partial v}{\partial t}
  = {}- U'(z - \xi(t))\, \xi'(t) - q'(t)
  = {}- V(U)\, \xi'(t) - q'(t) \,,
\\
\label{e:dv/dz}
&   \frac{\partial v}{\partial z}
  = U'(z - \xi(t)) = V(U) \,,
\\
\nonumber
&   f(v(z,t)) - f( U(z - \xi(t)) )
\\
\label{e:f(v)}
& = \Phi\left( U(z - \xi(t)) ,\, q(t) \right) \cdot q(t)
  = \Phi\left( U ,\, q(t) \right) \cdot q(t) \,,
  \quad\mbox{ where }\quad
\\
\label{e:Phi(v)}
& \Phi(s,q)\eqdef
  \frac{f(s-q) - f(s)}{q}
    \quad\mbox{ for }\, s\in \RR \,,\ q\in \RR\setminus \{ 0\} \,.
\end{align}
Consequently, the expressions in ineq.~\eqref{ineq:sub} become
\begin{align}
\label{eq:sub-sup}
&
\begin{aligned}
    \mathcal{N}(v) =
& {}- U'(z - \xi(t))\, \xi'(t) - q'(t)
    - U''(z - \xi(t)) - c\,  U'(z - \xi(t))
\\
& {}- f( U(z - \xi(t)) )
    - \Phi( U(z - \xi(t)) ,\, q(t) )\, q(t)
\end{aligned}
\\
\nonumber
& = {}- U'\, \xi' - q' - \Phi(U,q)\, q
  = {}- q\left[
    V(U)\, \frac{ \xi'(t) }{q(t)}
    + \frac{\mathrm{d}}{\mathrm{d}t}\, \ln q(t)
    + \Phi(U,q)
    \right] \,.
\end{align}

Notice that our choice of $q_0$ in \eqref{e:q_0} guarantees that
$0 < q_0\leq 1 - s_0 - 2\eta$.
Thus, we may apply ineq.~\eqref{ineq:f(1)} to conclude that
\begin{equation*}
  \Phi(s,q)\geq \overline{\mu}\geq \mu\eqdef
  \min\{ \underline{\mu} ,\, \overline{\mu} \} > 0
\end{equation*}
holds for all pairs $(s,q)$ satisfying
$1 - \delta\leq s\leq 1$ and $0 < q\leq q_0$.
Recall that the constant $\underline{\mu} > 0$
has been introduced in ineq.~\eqref{ineq:f(0)}.
Consequently, using eq.~\eqref{eq:sub-sup},
we will easily obtain ineq.~\eqref{ineq:sub}, that is,
$\mathcal{N}(v)\leq 0$ for all $(z,t)\in \Omega_{\delta}^{(3)}$,
provided we choose the functions
$\xi\colon \RR_+\to \RR$ and $q\colon \RR_+\to (0,\infty)$
as in \eqref{e:q_i,xi_i}, where
$\xi_{\infty,1} = \xi_{\infty}\in \RR$ and $\nu > 0$ are some constants,
such that
$\xi(0) = \xi_{\infty} - \nu q_0\geq z^{\ast}$
(to be specified later when we treat the second case).
From the derivatives
\begin{equation}
\label{e:q'/q,xi'/q}
    \frac{q'(t)}{q(t)}
  = \frac{\mathrm{d}}{\mathrm{d}t}\, \ln q(t)
  = {}- \mu \;(< 0) \quad\mbox{ and }\quad
    \frac{ \xi'(t) }{q(t)}
  = \mu\nu \;(> 0)
\end{equation}
inserted into the last bracket in eq.~\eqref{eq:sub-sup},
we deduce the desired ineq.~\eqref{ineq:sub}.

We continue with the \underline{first case}, $0\leq U\leq \delta$:
While treating the third case above ($1-\delta\leq U\leq 1$),
we have chosen $q_0\in \RR$ such that
inequalities \eqref{e:q_0} be satisfied.
In analogy with the third case, we wish to show that
$\mathcal{N}(v)\leq 0$ for all $(z,t)\in \Omega_{\delta}^{(1)}$.
We write the set
\begin{math}
  \Omega_{\delta}^{(1)} =
  \Omega_{\delta}^{(1,+)} \cup \Omega_{\delta}^{(1,-)}
\end{math}
as the union of the subsets
\begin{align}
\label{e:Omega_1+}
  \Omega_{\delta}^{(1,+)} &\eqdef
  \{ (z,t)\in \Omega_{\delta}^{(1)} \colon q(t)\leq U( z - \xi(t) ) \}
    \quad\mbox{ and }
\\
\label{e:Omega_1-}
  \Omega_{\delta}^{(1,-)} &\eqdef
  \left\{
    (z,t)\in \Omega_{\delta}^{(1)} \colon U( z - \xi(t) ) < q(t)
  \right\} \,.
\end{align}
The subsolution $v_1 = v$ in eq.~\eqref{e:v_1}
is now given by
\begin{equation}
\label{e:v_1:1}
  v(z,t) =
\left\{ \quad
\begin{alignedat}{2}
& U(z - \xi(t)) - q(t) &&\quad\mbox{ for }\,
    (z,t)\in \Omega_{\delta}^{(1,+)} \,,
\\
& 0 &&\quad\mbox{ for }\, 
    (z,t)\in \Omega_{\delta}^{(1,-)} \,.
\end{alignedat}
\right.
\end{equation}
Trivially, we have
$\mathcal{N}(0) = {}- f(0) = 0$.
It suffices to show that $\mathcal{N}(v)\leq 0$ holds for all
$(z,t)\in \Omega_{\delta}^{(1,+)}$.
This time we apply ineq.~\eqref{ineq:f(0)} with
$\delta\in (0,\eta)$ and $\underline{\mu} > 0$ to conclude that
\begin{equation*}
  \Phi(s,q)\geq \underline{\mu}\geq \mu =
  \min\{ \underline{\mu} ,\, \overline{\mu} \} > 0
\end{equation*}
holds for all pairs $(s,q)$ satisfying
$0 < q\leq s\leq \delta$.
We recall the inequalities
$0 < \delta < s_0 - 2\eta$, by \eqref{e:delta<1/3} and
$0 < \delta < \eta\leq \eta_0$.
Consequently, using eq.~\eqref{eq:sub-sup}, we finally obtain
$\mathcal{N}(v)\leq 0$ for all $(z,t)\in \Omega_{\delta}^{(1,+)}$
in the same way as in the third case above
($1-\delta\leq U\leq 1$).

It remains to treat the \underline{second case},
$\delta\leq U\leq 1 - \delta$:
Here, we take advantage of the fact that there is a constant
$\omega\in (0,\infty)$ such that $V(s)\geq \omega$ for all
$s\in [\delta, 1-\delta]$, by ineq.~\eqref{e:dU/dz}.
Furthermore, ineq.~\eqref{e:f_Lip}
(a one\--sided Lipschitz condition)
guarantees that $\Phi(s,q)\geq {}- L$ whenever
$0 < q\leq s\leq 1$.
Similarly to the first case above ($0\leq U\leq \delta$),
we write the set
\begin{math}
  \Omega_{\delta}^{(2)} =
  \Omega_{\delta}^{(2,+)} \cup \Omega_{\delta}^{(2,-)}
\end{math}
as the union of the subsets
\begin{align}
\label{e:Omega_2+}
  \Omega_{\delta}^{(2,+)} \eqdef
& \left\{
    (z,t)\in \Omega_{\delta}^{(2)} \colon q(t)\leq U( z - \xi(t) )
  \right\}
    \quad\mbox{ and }
\\
\label{e:Omega_2-}
  \Omega_{\delta}^{(2,-)} \eqdef
& \left\{
    (z,t)\in \Omega_{\delta}^{(2)} \colon U( z - \xi(t) ) < q(t)
  \right\} \,.
\end{align}
The subsolution $v_1 = v$ in eq.~\eqref{e:v_1}
is given by
\begin{equation}
\label{e:v_1:2}
  v(z,t) =
\left\{ \quad
\begin{alignedat}{2}
& U(z - \xi(t)) - q(t) &&\quad\mbox{ for }\,
    (z,t)\in \Omega_{\delta}^{(2,+)} \,,
\\
& 0 &&\quad\mbox{ for }\, 
    (z,t)\in \Omega_{\delta}^{(2,-)} \,.
\end{alignedat}
\right.
\end{equation}
Again, as in the first case, it suffices to show that
$\mathcal{N}(v)\leq 0$ holds for all
$(z,t)\in \Omega_{\delta}^{(2,+)}$.
This time we apply eq.~\eqref{e:q'/q,xi'/q}
and the inequalities
\begin{alignat*}{2}
  V(s) &\geq \omega\equiv \mathrm{const} > 0
&&  \quad\mbox{ for all }\, s\in [\delta, 1-\delta] \,,
\\
  \Phi(s,q) &\geq {}- L\equiv \mathrm{const}\leq 0
&&  \quad\mbox{ whenever }\, 0 < q\leq s\leq 1 \,,
\end{alignat*}
obtained above to eq.~\eqref{eq:sub-sup}, thus arriving at
\begin{align*}
    \mathcal{N}(v)
& = {}- q\left[
    V(U)\, \frac{ \xi'(t) }{q(t)}
    + \frac{\mathrm{d}}{\mathrm{d}t}\, \ln q(t)
    + \Phi(U,q)
    \right]
\\
& \leq {}- q\left[ \omega\, \mu\nu - \mu - L\right]
  \leq 0
  \quad\mbox{ for all }\, (z,t)\in \Omega_{\delta}^{(2,+)} \,,
\end{align*}
provided $\nu > 0$ is chosen sufficiently large, such that
$[\,\cdots\,] = \omega\, \mu\nu - \mu - L\geq 0$.
For instance,
\begin{equation*}
  \nu = \frac{ 1 + (L/\mu) }{\omega} > 0
\end{equation*}
will do it.
Moreover, the remaining unspecified constant
$\xi_{\infty}$ from the third case above may be chosen to be
\begin{equation*}
    \xi_{\infty} = \nu q_0 + z^{\ast}
  = \frac{ 1 + (L/\mu) }{\omega}\cdot q_0
  + z^{\ast} \in \RR \,.
\end{equation*}
We have verified that $v_1 = v$ given by eq.~\eqref{e:v_1:2}
is a subsolution to the Cauchy problem \eqref{e_c:FKPP}
also in the second case
($\delta\leq U\leq 1 - \delta$).

This proves that the function $v_1 = v$ defined in eq.~\eqref{e:v_1}
is a subsolution to the Cauchy problem \eqref{e_c:FKPP}
in all three cases above.
Moreover, $v$ is Lipschitz\--continuous in each set
$\Omega_{\delta}^{(j)}$; $j=1,2,3$, and, consequently,
by its definition in \eqref{e:subsol},
also in the whole space\--time domain
\begin{math}
  \RR\times \RR_+ =
  \Omega_{\delta}^{(1)}\cup \Omega_{\delta}^{(2)}\cup
  \Omega_{\delta}^{(3)} \,.
\end{math}
It follows that $v$ is a weak $L^{\infty}$-subsolution to
\eqref{e_c:FKPP}.

A weak $L^{\infty}$-supersolution $v_2$ to
the Cauchy problem \eqref{e_c:FKPP}
is obtained analogously, cf.\ eq.~\eqref{e:supsol};
we leave the details to an interested reader.
We remark that the constant
\begin{math}
  \mu = \min\{ \underline{\mu} ,\, \overline{\mu} \} > 0
\end{math}
remains the same for both $v_i$; $i=1,2$.
\qed
\par\vskip 10pt

\section{Long\--Time Asymptotic Behavior of Solutions}
\label{s:Time-Asympt}

This section is concerned with
the long\--time asymptotic behavior of solutions $v = v(z,t)$
to the initial value problem \eqref{e_c:FKPP} as time $t\to \infty$.
More precisely, the uniform convergence of the family of functions
$v(\,\cdot\,,t)\colon \RR\to [0,1]$; $t\in \RR_+$,
to a TW\--solution
$z\mapsto U(z + \zeta)\colon \RR\to [0,1]$ as $t\to \infty$
will be proved in our main result, Theorem~\ref{thm-Lyapunov}.
Our approach is based on controlling
the long\--time asymptotic behavior of $v(z,t)$ pointwise by
a special pair of sub- and super\-solutions constructed in
Proposition~\ref{prop-sub-sup}.
We combine the Lyapunov stability of TW\--solutions with
the spatial regularity of solutions $v(z,t)$ and
the minimization of the Lyapunov functional
on the $\omega$\--limit set of a solution as time goes to infinity
in order to prove the long\--time convergence
in Theorem~\ref{thm-Lyapunov}.

\subsection{Lyapunov Stability of Travelling Waves}
\label{ss:Stability}

In this paragraph we establish the stability of travelling waves
(equivalently, that of TW\--solutions)
in the sense of Lyapunov.
These travelling waves for problem \eqref{e:FKPP},
$(x,t)\mapsto u(x,t) = U(x-ct + \zeta)$,
where $\zeta\in \RR$ is arbitrary,
have been obtained in Section~\ref{s:Trav_Waves}; see
Proposition~\ref{prop-uniq_TW}.
Roughly speaking, if the (bounded classical) solution of
problem~\eqref{e_c:FKPP}, say,
$v(\,\cdot\,,t)\colon \RR\to [0,1]$; $t\in \RR_+$,
has approached a TW\--solution
$z\mapsto U(z + \zeta)$ ``close enough''
at some time $t_0\in \RR_+$ and for some $\zeta\in \RR$,
that is to say, if the difference\hfil\break
$v(z,t_0) - U(z + \zeta)$ is ``small enough'',
uniformly for all $z\in \RR$, then the difference\hfil\break
$v(z,t) - U(z + \zeta)$ stays ``small''
for all later times $t\geq t_0$, uniformly for all $z\in \RR$.
This claim, called {\em Lyapunov stability\/}, is quantified below.

\begin{proposition}\label{prop-Lyapunov}
Assume that\/ $f$ satisfies all four\/ {\rm Hypotheses},
{\bf (H1)} through {\bf (H4)}.
Let\/
$z\mapsto U(z + \zeta)\colon \RR\to [0,1]$,
$\zeta\in \RR$, be a TW\--solution as specified in
{\rm Proposition~\ref{prop-uniq_TW}, Part~{\rm (d)}}.
Then there is a function
$\varrho = \varrho(\eps) > 0$ defined for every\/
$\eps\in (0,\eps_0]$, $\eps_0 > 0$ small enough,
\begin{math}
  \eps_0\leq \frac{1}{2}\cdot \min\{ s_0 ,\, 1 - s_0\} \,,
\end{math}
such that\/
$\lim_{\eps\to 0+} \varrho(\eps) = 0$ and the following property holds:
\begin{itemize}
\item[{\rm (LS)}]
Let\/ $\eps\in (0,\eps_0)$ be arbitrary.
If\/ $v_0\in L^{\infty}(\RR)$ satisfies\/
\begin{equation}
\label{e:v_0-U}
  0\leq v_0(z)\leq 1 \quad\mbox{ and }\quad
  |v_0(z) - U(z + \zeta)| \leq \eps
    \quad\mbox{ for all }\, z\in \RR \,,
\end{equation}
where $\zeta\in \RR$ is a suitable constant,
then the (unique bounded classical) solution
$v\colon \RR\times \RR_+\to \RR$ of problem \eqref{e_c:FKPP}
{\rm (cf.\ Proposition~\ref{prop-weak_sol})}
satisfies\/
\begin{equation}
\label{e:v(t)-U}
  |v(z,t) - U(z + \zeta)| \leq \varrho(\eps)
    \quad\mbox{ for all }\, (z,t)\in \RR\times \RR_+ \,.
\end{equation}
\end{itemize}
\end{proposition}
\par\vskip 10pt

Since $v_0\in L^{\infty}(\RR)$ in {\rm Condition (LS)} satisfies
inqualities \eqref{e:v_0-U} with
\begin{math}
  0 < \eps\leq \eps_0\leq \frac{1}{2}\cdot \min\{ s_0 ,\, 1 - s_0\} \,,
\end{math}
it obeys also {\rm Hypothesis {\bf (H5)}}, with a help from
\eqref{lim:FKPP}.

\par\vskip 10pt
\proof
Set
\begin{math}
  \eps_0' = \frac{1}{2}\cdot \min\{ s_0 ,\, 1 - s_0\} > 0 \,.
\end{math}
We begin with an arbitrary number $\eps\in (0, \eps_0']$.
Let $v_0\in L^{\infty}(\RR)$ satisfy inequalities \eqref{e:v_0-U}
with some constant $\zeta\in \RR$.
Next, we recall that, by Proposition~\ref{prop-sub-sup},
there is an orderded pair of
weak $L^{\infty}$-sub- and -super\-solutions of
the Cauchy problem \eqref{e_c:FKPP},
$v_1(z,t)$ and $v_2(z,t)$, respectively, given by formulas
\eqref{e:subsol} and \eqref{e:supsol}, such that
$0\leq v_1(z,t)\leq v_2(z,t)\leq 1$ for all
$(z,t)\in \RR\times \RR_+$, and
$v_1(z,0)\leq v_0(z)\leq v_2(z,0)$ for all $z\in \RR$.
With regard to inequalities \eqref{e:v_0-U},
we are going to find the constants
$\mu, \nu, q_{0,i}\in (0,\infty)$ and
$\xi_{\infty,i}\in \RR$ in Proposition~\ref{prop-sub-sup}
($i=1,2$), such that we have also
\begin{equation}
\label{e:v_1<U-eps<U+eps<v_2}
\begin{aligned}
  v_1(z,0)
&         \leq \max\{ U(z + \zeta) - \eps ,\, 0\}
          \leq v_0(z)
\\
&         \leq \min\{ U(z + \zeta) + \eps ,\, 1\}
          \leq v_2(z,0)
  \quad\mbox{ for }\, z\in \RR
\end{aligned}
\end{equation}
(at $t=0$).
Recalling eqs.~\eqref{e:q_i,xi_i} we observe that inequalities
\eqref{e:v_1<U-eps<U+eps<v_2} above are satisfied provided
the following two inequalities hold for every $z\in \RR$:
\begin{alignat}{2}
\label{e:v_1<U-eps}
  U(z - \xi_{\infty,1}) - q_{0,1}
& \leq \max\{ U(z + \zeta) - \eps ,\, 0\}
  \quad\mbox{ and }\quad
\\
\label{e:U+eps<v_2}
  \min\{ U(z + \zeta) + \eps ,\, 1\}
& \leq U(z + \xi_{\infty,2}) + q_{0,2} \,.
\end{alignat}
Clearly, by our choice of $\eps\in (0, \eps_0']$,
both functions of $z\in \RR$,
\begin{equation}
\label{e:v_0,1:v_0,2}
  v_{0,1}(z) = \max\{ U(z + \zeta) - \eps ,\, 0\}
  \quad\mbox{ and }\quad
  v_{0,2}(z) = \min\{ U(z + \zeta) + \eps ,\, 1\} \,,
\end{equation}
satisfy {\rm Hypothesis {\bf (H5)}} for the initial condition $v_0$,
in addition to $v_{0,1}\leq v_0\leq v_{0,2}$ on $\RR$.
Consequently, we may apply Proposition~\ref{prop-sub-sup}
to find some constants
$\mu, \nu, q_{0,i}\in (0,\infty)$ and
$\xi_{\infty,i}\in \RR$ ($i=1,2$), such that both inequalities
\eqref{e:v_1<U-eps} and \eqref{e:U+eps<v_2} are valid, provided
$\eps$ is sufficiently small, say, $0 < \eps\leq \eps_0'' < \infty$.
Setting
\begin{math}
  \eps_0 = \min\{ \eps_0' ,\, \eps_0''\} > 0
\end{math}
and inspecting our proof of Proposition~\ref{prop-sub-sup},
we conclude that the constants $\mu, \nu\in (0,\infty)$
can be chosen independent from $\eps\in (0,\eps_0]$;
$\mu > 0$ sufficiently small and $\nu > 0$ sufficiently large.
In contrast, $q_{0,i}\in (0,\infty)$ and $\xi_{\infty,i}\in \RR$
depend on $\eps\in (0,\eps_0]$ in such a way that the estimates
\begin{align}
\label{e:q_0,1:q_0,2}
& 0 < q_{0,i}\equiv q_{0,i}(\eps)\leq C\eps \,;\quad i=1,2 \,,
  \quad\mbox{ and }\quad
\\
\label{e:xi_0,1:xi_0,2}
& |\xi_{\infty,1} + \zeta| \leq C\eps \,,\quad
  |\xi_{\infty,2} - \zeta| \leq C\eps
\end{align}
hold with a constant $C > 0$ independent from $\eps\in (0,\eps_0]$.
From the formulas in eqs.~\eqref{e:q_i,xi_i} we derive also
\begin{align}
\label{e:q_1:q_2}
& 0 < q_i(t)\leq q_{0,i}\leq C\eps \,;\quad i=1,2 \,,
\\
\label{e:xi_1(t)}
& |\xi_1(t) + \zeta| \leq
  |\xi_{\infty,1} + \zeta| + \nu\, q_1(t)
  \leq (1 + \nu) C\eps \,,\quad\mbox{ and }\quad
\\
\label{e:xi_2(t)}
& |\xi_2(t) - \zeta| \leq
  |\xi_{\infty,2} - \zeta| + \nu\, q_2(t)
  \leq (1 + \nu) C\eps
\end{align}
for all $t\in \RR_+$ and for every $\eps\in (0,\eps_0]$.
We combine these inequalities with the derivative
$U'\colon \RR\to \RR$ being bounded to conclude that,
by formulas \eqref{e:subsol} and \eqref{e:supsol},
there is a constant $C'> 0$ independent from $\eps\in (0,\eps_0]$,
such that
\begin{equation}
\label{e:v_1-v_2}
  0\leq v_2(z,t) - v_1(z,t)\leq C'\eps
    \quad\mbox{ holds for all }\, (z,t)\in \RR\times \RR_+ \,.
\end{equation}
Since also all inequalities
$v_1(z,t)\leq v(z,t)\leq v_2(z,t)$ and
$v_1(z,t)\leq U(z + \zeta)\leq v_2(z,t)$ are valid for all
$(z,t)\in \RR\times \RR_+$,
by Lemma~\ref{lem-weak_comp} (weak comparison principle),
we arrive at the desired estimate \eqref{e:v(t)-U}
with the function $\varrho(\eps) = C'\eps$ defined for every
$\eps\in (0,\eps_0]$.
\qed
\par\vskip 10pt

\subsection{Compactness and Regularity of Solution Orbits}
\label{ss:comp_orbits}

The regularity and relative compactness of orbits generated by
the unique solutions of the initial value problem \eqref{e_c:FKPP}
are treated in the next two lemmas.

\begin{lemma}\label{lem-time_reg}
Let\/ $T\in (0,\infty)$ and let\/
$\alpha\in (0,1)$ be the H\"older exponent for\/ $f$ from
{\rm Hypothesis {\bf (H2)}}.
Given any Lebesgue measurable initial data
$v_0\colon \RR\to [0,1]$, let\/
$v\colon \RR\times (0,\infty)\to \RR$
denote the unique bounded classical solution of problem~\eqref{e_c:FKPP}.
Let us denote by
\begin{math}
  v_{\tau}(z,t)\eqdef v(z, t + \tau) \,,
    \;\mbox{ for }\, (z,t)\in \RR\times [0,T] \,,
\end{math}
the time translation of $v$ by $\tau$ time units,
$\tau\in (0,\infty)$.
Then, for any $\tau_0\in (0,\infty)$ there exists a constant\/
$C_{(0,T)}^{\alpha}\equiv C_{(0,T)}^{\alpha}(\tau_0)\in \RR_+$
such that the estimate\/
\begin{align}
\label{tau:Hoelder:v_2}
& \max\left\{
  \left\Vert \frac{\partial v_{\tau}}{\partial z}
  \right\Vert_{ C^{\alpha, \alpha/2}(\RR\times [0,T]) } \,,\;
  \left\Vert \frac{\partial^2 v_{\tau}}{\partial z^2}
  \right\Vert_{ C^{\alpha, \alpha/2}(\RR\times [0,T]) } \,,\;
  \left\Vert \frac{\partial v_{\tau}}{\partial t}
  \right\Vert_{ C^{\alpha, \alpha/2}(\RR\times [0,T]) }
    \right\}
\\
\nonumber
& \leq C_{(0,T)}^{\alpha}
\end{align}
holds for every time translation $\tau\geq \tau_0 > 0$.
\end{lemma}

\par\vskip 10pt
\proof
A proof of ineq.~\eqref{tau:Hoelder:v_2}
follows directly from the regularity results in \eqref{Hoelder:v_2}.
\qed
\par\vskip 10pt

We combine Lemma~\ref{lem-time_reg} with
Arzel\`a\--Ascoli's compacness criterion to obtain the following result.

\begin{lemma}\label{lem-time_comp}
Let\/ $T\in (0,\infty)$ and let\/
$\alpha\in (0,1)$ be the H\"older exponent for\/ $f$ from
{\rm Hypothesis {\bf (H2)}}.
Assume that\/ $J\subset \RR$ is a compact interval, $\beta\in (0,\alpha)$,
and\/ $0 < \tau_0 < \infty$.
Under the hypotheses of\/ {\rm Lemma~\ref{lem-time_reg}},
the family of ``translation'' functions\/
$v_{\tau}\colon J\times [0,T]\to \RR$; $\tau\geq \tau_0$,
is relatively compact in the H\"older space
$C^{2+\beta, 1+(\beta/2)}(J\times [0,T])$.
\end{lemma}
\par\vskip 10pt

\subsection{Convergence of Solutions to a Travelling Wave}
\label{ss:Trav-Asympt}

From now on we treat solely the case
$-\infty < z_0 < z_1 < \infty$
as opposed to the classical case
$z_0 = -\infty$ and $z_1 = +\infty$
($f\colon \RR\to \RR$ continuously differentiable)
treated in
{\sc P.~C.\ Fife} and {\sc J.~B.\ Mc{L}eod}
\cite[Sect.~1]{Fife-McLeod} and
{\sc J.~D.\ Murray} \cite{Murray-I}, {\S}13.3, pp.\ 444--449.
In {\S}\ref{ss:Trav-Shape} we have given a simple example when
$-\infty < z_0 < z_1 < \infty$ (Example~\ref{exam-z_0<z_1}).

Let us recall that we work with initial data
$v_0\colon \RR\to \RR$ satisfying
{\rm Hypothesis {\bf (H5)}}.
Given these initial data, in Proposition~\ref{prop-sub-sup}
we have constructed an ordered pair of
weak $L^{\infty}$-sub- and -supersolutions,
$v_1$ and $v_2$, respectively, such that
$v_1(z,0)\leq v_0(z)\leq v_2(z,0)$ for all $z\in \RR$ at $t=0$, and
$0\leq v_1(z,t)\leq v_2(z,t)\leq 1$ for all
$(z,t)\in \RR\times \RR_+$.
Especially formulas in \eqref{e:q_i,xi_i} are of importance.
Recalling our definition of $z_0, z_1$ in \eqref{lim:dz/dU},
we choose a compact interval
$J = [a,b]\subset \RR$ such that
\begin{equation}
\label{ineq:a+xi_2<z_0}
  a + \xi_{\infty,2} < z_0 < z_1 < b - \xi_{\infty,1} \,.
\end{equation}
Since $\xi_i(t)\leq \xi_{\infty,i}$ for $t\geq 0$; $i=1,2$,
this choice of $J$ guarantees
\begin{equation}
\label{ineq:a+xi_2(t)<z_0}
  a + \xi_2(t) < z_0 < z_1 < b - \xi_1(t)
    \quad\mbox{ for all }\, t\geq 0 \,,
\end{equation}
together with the limits
\begin{align}
\label{lim:a+xi_2(t)<z_0}
  \lim_{t\to \infty} v_1(z,t) =
& \lim_{t\to \infty} v_2(z,t) = 0
  \quad\mbox{ uniformly for }\, z\in (-\infty,a] \,,
  \quad\mbox{ and }
\\
\label{lim:z_1<b-xi_1(t)}
  \lim_{t\to \infty} v_1(z,t) =
& \lim_{t\to \infty} v_2(z,t) = 1
  \quad\mbox{ uniformly for }\, z\in [b,+\infty) \,.
\end{align}
%

\begin{lemma}\label{lem-time_conv}
Let\/ $T\in (0,\infty)$.
Assume that\/ $f$ satisfies all four\/ {\rm Hypotheses},
{\bf (H1)} through {\bf (H4)}.
Let\/ $J\subset \RR$ be the compact interval specified above,
$\beta\in (0,\alpha)$, and\/ $0 < \tau_0 < \infty$.
Assume that the initial data $v_0\colon \RR\to \RR$ satisfy\/
{\rm Hypothesis {\bf (H5)}} and let\/
$v\colon \RR\times (0,\infty)\to \RR$
denote the unique bounded classical solution of problem~\eqref{e_c:FKPP}.
Finally, let\/
$v_{\tau}\colon J\times [0,T]\to \RR$
be the family of ``translation'' functions defined in
{\rm Lemma~\ref{lem-time_reg}} for $\tau\geq \tau_0$.
Then every sequence
$\{\tau_n\}_{n=1}^{\infty} \subset [\tau_0,\infty)$,
$\tau_n\to \infty$ as $n\to \infty$,
contains a subsequence, denoted again by
$\{\tau_n\}_{n=1}^{\infty}$, such that\/
\begin{equation}
\label{v_tau_n->v*}
  v_{\tau_n}\to v^{*} \quad\mbox{ in }\,
  C^{2+\beta, 1+(\beta/2)}(J\times [0,T])
    \quad\mbox{ as }\, n\to \infty \,,
\end{equation}
for some function
$v^{*} \in C^{2+\beta, 1+(\beta/2)}(J\times [0,T])$.

\begin{enumerate}
\renewcommand{\labelenumi}{{\bf (H\alph{enumi})}}
\item[{\rm (i)}]
In particular, we have $v^{*}(a,t) = 0$, $v^{*}(b,t) = 1$, and\/
\begin{equation}
\label{e:U_1<v*<U_2}
  U(z - \xi_{\infty,1})\leq v^{*}(z,t)\leq U(z + \xi_{\infty,2})
  \quad\mbox{ for all }\, (z,t)\in J\times [0,T] \,.
\end{equation}
\item[{\rm (ii)}]
The limit function $v^{*}\colon J\times [0,T]\to \RR$
satisfies eq.\ \eqref{e_c:FKPP} in $J\times [0,T]$ in the classical sense
(with all partial derivatives being continuous).
Consequently, its natural extension $\tilde{v}$ to $\RR\times [0,T]$,
\begin{equation*}
  \tilde{v}(z,t) =
\left\{
  \begin{array}{cl}
  v^{*}(z,t) &\quad\mbox{ if }\, z\in J = [a,b] \,;
\\
  0 &\quad\mbox{ if }\, z\in (-\infty,a) \,;
\\
  1 &\quad\mbox{ if }\, z\in (b,\infty) \,,
  \end{array}
\right.
\end{equation*}
defined for all $(z,t)\in \RR\times [0,T]$,
is a classical solution to eq.\ \eqref{e_c:FKPP} in $\RR\times [0,T]$.
\item[{\rm (iii)}]
We have the uniform convergence on all of\/ $\RR\times [0,T]$,
\begin{equation}
\label{lim:v_tau_n->v*}
  \sup_{ (z,t)\in \RR\times [0,T] }\
  \left| v(z, t + \tau_n) - \tilde{v}(z,t) \right|
    \,\longrightarrow\, 0 \quad\mbox{ as }\, n\to \infty \,.
\end{equation}
\end{enumerate}
\end{lemma}

\par\vskip 10pt
\proof
The convergence result in \eqref{v_tau_n->v*}
follows directly from Lemma~\ref{lem-time_comp} above.

Proof of {\rm (i)}.
The inequalities in \eqref{e:U_1<v*<U_2} are derived from
the initial data $v_0$ satisfying
$v_1(z,0)\leq v_0(z)\leq v_2(z,0)$ for all $z\in \RR$ at $t=0$,
by {\rm Hypothesis {\bf (H5)}} combined with
Proposition~\ref{prop-sub-sup} and
the weak comparison principle in Lemma~\ref{lem-weak_comp};
first, for each function $v_{\tau_n}$; $n=1,2,3,\dots$,
\begin{equation*}
  v_1(z, t + \tau_n)\leq v_{\tau_n}(z,t)\leq v_2(z, t + \tau_n)
  \quad\mbox{ for all }\, (z,t)\in J\times [0,T] \,,
\end{equation*}
then for the limit function $v^{*}$ as $n\to \infty$.

The boundary behavior of $v^{*}$, i.e.,
$v^{*}(a,t) = 0$, $v^{*}(b,t) = 1$,
is obtained from \eqref{ineq:a+xi_2(t)<z_0} and \eqref{e:U_1<v*<U_2}.
From this boundary behavior of $v^{*}$,
\eqref{v_tau_n->v*}, and \eqref{e:U_1<v*<U_2} we derive also
\begin{equation*}
    \frac{\partial v^{*}}{\partial z} (a,t)
  = \frac{\partial v^{*}}{\partial z} (b,t) = 0
    \quad\mbox{ and }\quad
    \frac{\partial^2 v^{*}}{\partial z^2} (a,t)
  = \frac{\partial^2 v^{*}}{\partial z^2} (b,t) = 0 \,.
\end{equation*}

Proof of {\rm (ii)}.
Each function $v_{\tau_n}$ is a bounded classical solution of
the equation in \eqref{e_c:FKPP}.
The convergence result in \eqref{v_tau_n->v*} yields that so is
the limit function $v^{*}$.
The remaining claims for $\tilde{v}$ follow from~{\rm (i)} and its proof.

Proof of {\rm (iii)}.
The uniform convergence on $J\times [0,T]$ follows from
\eqref{v_tau_n->v*}.
To verify it also on the complement
\begin{equation*}
    (\RR\setminus J)\times [0,T]
  = \bigl( (-\infty,a)\times [0,T] \bigr) \cup
    \bigl( (b,+\infty)\times [0,T] \bigr) ,
\end{equation*}
let us recall our choice of the compact interval
$J = [a,b]\subset \RR$ such that \eqref{ineq:a+xi_2<z_0} is valid.
Then formulas \eqref{e:subsol} and \eqref{e:supsol} yield
$v_1(z,0)\leq v_0(z)\leq v_2(z,0)$ for all $z\in \RR$, and
$v_1(z,t)\leq v(z,t)\leq v_2(z,t)$ for all
$(z,t)\in \RR\times \RR_+$.
Recalling\hfil\break
\begin{math}
  v_{\tau_n}(z,t)\eqdef v(z, t + \tau_n)
    \quad\mbox{ for }\, (z,t)\in \RR\times [0,T] \,,
\end{math}
and letting $n\to \infty$, we have also
\begin{equation}
\label{lim:v(tau_n)}
\begin{aligned}
&      U(z - \xi_{\infty,1}) = \lim_{n\to \infty} v_1(z, t + \tau_n)
  \leq \tilde{v}(z,t)
\\
& \leq \lim_{n\to \infty} v_2(z, t + \tau_n) = U(z + \xi_{\infty,2})
    \quad\mbox{ for }\, (z,t)\in (\RR\setminus J)\times [0,T] \,.
\end{aligned}
\end{equation}
We combine these estimates with
$v_1\leq v\leq v_2$ in $\RR\times \RR_+$, thus arriving at
\begin{align*}
& \lim_{n\to \infty}\
  \sup_{ (z,t)\in (\RR\setminus J)\times [0,T] }\
  \left| v(z, t + \tau_n) - \tilde{v}(z,t) \right|
\\
& \leq \sup_{ z\in \RR\setminus J }\
    \bigl( U(z + \xi_{\infty,2}) - U(z - \xi_{\infty,1}) \bigr)
  = 0 \,,
\end{align*}
thanks to \eqref{ineq:a+xi_2<z_0}.
This proves the uniform convergence on the complement
$(\RR\setminus J)\times [0,T]$.
We combine it with the uniform convergence on $J\times [0,T]$
in \eqref{v_tau_n->v*} to derive \eqref{lim:v_tau_n->v*}
in Part~{\rm (iii)}.
\qed
\par\vskip 10pt

In our treatment of the dynamical system generated by the solutions of
problem \eqref{e_c:FKPP}
we use standard terminology from {\sc J.~K.\ Hale} \cite{Hale}.
The trivialized analogue of Lemma~\ref{lem-time_conv} for $T=0$
enables us to define the {\em $\omega$\--limit set\/}
$\omega(v_0)$ of the solution $v$ as follows:
Given any $0 < \tau < \infty$, the set of functions
\begin{equation*}
  \mathcal{O}_{\tau} =  \{ v(\,\cdot\,,t)\colon t\geq \tau\}
\end{equation*}
has a compact closure
$\overline{\mathcal{O}_{\tau}}$ in $C^2(J)$.
We define
\begin{equation*}
  \omega(v_0)\eqdef \bigcap_{\tau > 0} \overline{\mathcal{O}_{\tau}} \,,
\end{equation*}
which is a nonempty compact set in $C^2(J)$.
It is important to notice that all conclusions of
Lemma~\ref{lem-time_conv}, {\rm Part~(i)},
remain valid also for every function $w\in \omega(v_0)$, whence
\begin{equation}
\label{e:U_1<w<U_2}
  w(a) = 0 \,,\ w(b) = 1 \,,\ \mbox{ and }\
    \frac{\mathrm{d}w}{\mathrm{d}z}(a)
  = \frac{\mathrm{d}w}{\mathrm{d}z}(b) = 0 \,.
\end{equation}
It is well\--known that the $\omega$\--limit set $\omega(v_0)$
is invariant under the mapping
$v^{*}(\,\cdot\,,0)\mapsto v^{*}(\,\cdot\,,t)$
for any $t\in [0,T]$; see~\cite{Hale}.
Next, we find an $\omega$\--limit point
$v^{*}(\,\cdot\,,0)\in \omega(v_0)$ such that
$v^{*}(\,\cdot\,,t)\equiv v^{*}(\,\cdot\,,0)$
in $J$ for all $t\in [0,T]$.

We define a Lyapunov functional on $\omega(v_0)$ as follows:
\begin{equation}
\label{def:Lyapunov}
  \mathcal{E}(w)\eqdef \int_a^b
  \biggl[ \frac{1}{2} \genfrac{(}{)}{}0{\mathrm{d}w}{\mathrm{d}z}^2
       - F(w(z))
  \biggr] \,\ee^{cz} \,\mathrm{d}z
  \quad\mbox{ for }\, w\in \omega(v_0) \,.
\end{equation}
This functional is continuous on the nonempty compact set
$\omega(v_0)$ in $C^2(J)$.
Combining the equation in \eqref{e_c:FKPP} for $v^{*}$
with the boundary behavior \eqref{e:U_1<w<U_2} of
$v^{*}(\,\cdot\,,t)\in \omega(v_0)$, we arrive at
\begin{equation}
\label{e:d_Lyapunov/dt}
  \frac{\mathrm{d}}{\mathrm{d}t}\,
    \mathcal{E}(v^{*}(\,\cdot\,,t))
  = {}- \int_a^b \genfrac{(}{)}{}0{\partial v^{*}}{\partial t}^2
    \,\ee^{cz} \,\mathrm{d}z
  \leq 0 \quad\mbox{ for }\, t\in [0,T] \,,
\end{equation}
cf.\ 
{\sc P.~C.\ Fife} and {\sc J.~B.\ Mc{L}eod}
\cite[p.~350]{Fife-McLeod}, proof of Lemma 4.5, or
{\sc J.~K.\ Hale} \cite[p.~76]{Hale}, {\S}4.3.1.
If $w_0\in \omega(v_0)$ is a minimizer for
$\mathcal{E}\colon \omega(v_0)\to \RR$
and we choose the sequence
$\{\tau_n\}_{n=1}^{\infty} \subset [\tau_0,\infty)$
in Lemma~\ref{lem-time_conv} such that
$v(\,\cdot\,, \tau_n)\to w_0$ in $C^2(J)$ as $n\to \infty$,
then we have also
\begin{equation}
\label{C^2,1:v_tau_n->v*}
  v_{\tau_n}\to v^{*} \quad\mbox{ in }\,
  C^{2,1}(J\times [0,T])
    \quad\mbox{ as }\, n\to \infty \,,
\end{equation}
by \eqref{v_tau_n->v*} in Lemma~\ref{lem-time_conv}, and
$v^{*}(\,\cdot\,,0) = w_0\in C^2(J)$.
If the function
\begin{math}
  t\mapsto v^{*}(\,\cdot\,,t)\colon [0,T]\to \omega(v_0)
  \subset C^2(J)
\end{math}
were not constant in time $t$, then
eq.~\eqref{e:d_Lyapunov/dt} would yield
$\mathcal{E}(v^{*}(\,\cdot\,,T)) < \mathcal{E}(w_0)$.
This contradicts our choice of $w_0\in \omega(v_0)$
to be a minimizer for
$\mathcal{E}\colon \omega(v_0)\to \RR$.
Hence, we have
$v^{*}(\,\cdot\,,t)\equiv w_0$ in $J$ for all $t\in [0,T]$.
The limit function $v^{*}(z,t)\equiv w_0(z)$
satisfies the equation in \eqref{e_c:FKPP} for $v^{*}$ with
${\partial v^{*}}/{\partial t} \equiv 0$ in $J\times [0,T]$, that is,
eq.~\eqref{eq:v(z)} together with \eqref{e:z<0,z>1},
by {\rm Proposition~\ref{prop-monot_TW}},
where the open interval $(z_0,z_1)$
has to be replaced by another open interval
$(\tilde{z}_0, \tilde{z}_1)\subset \RR$.
We apply {\rm Proposition~\ref{prop-uniq_TW}}, Part~{\rm (d)},
to specify the shape of $w_0$:
There is a number $\zeta\in \RR$ such that
$w_0(z) = U(z + \zeta)$ for all $z\in \RR$, and
\begin{math}
  (\tilde{z}_0, \tilde{z}_1) = (z_0 - \zeta, z_1 - \zeta)
  \subset J = [a,b] .
\end{math}
Consequently, \eqref{C^2,1:v_tau_n->v*} reads
\begin{equation}
\label{z,t:v_tau_n->v*}
  \left\| v(z, t + \tau_n) - U(z + \zeta)
  \right\|_{ C^{2,1}(J\times [0,T]) } \;\longrightarrow\; 0
    \quad\mbox{ as }\, n\to \infty \,,
\end{equation}
where the norm is taken for the function
\begin{math}
  (z,t)\mapsto v(z, t + \tau_n) - U(z + \zeta)
  \colon J\times [0,T]\to \RR .
\end{math}
We would like to emphasize that the shift
$\zeta\equiv \zeta (v_0)$
depends on the choice of the initial data $v_0$
and on our choice of the $\omega$\--limit point $w_0\in \omega(v_0)$.
In turn, these choices specify also the sequence (or subsequence of)
$\{\tau_n\}_{n=1}^{\infty} \subset [\tau_0,\infty)$.

We collect the results of this section up to now in
the following proposition:

\begin{proposition}\label{prop-time_conv}
Assume that\/ $f$ satisfies all four\/ {\rm Hypotheses},
{\bf (H1)} through {\bf (H4)}.
Let\/ $J\subset \RR$ be the compact interval specified above,
and\/ $0 < \tau_0 < \infty$.
Assume that the initial data $v_0\colon \RR\to \RR$ satisfy\/
{\rm Hypothesis {\bf (H5)}} and let\/
$v\colon \RR\times (0,\infty)\to \RR$
denote the unique bounded classical solution of problem~\eqref{e_c:FKPP}.
Then there exists a sequence
$\{\tau_n\}_{n=1}^{\infty} \subset [\tau_0,\infty)$,
$\tau_n\to \infty$ as $n\to \infty$
such that, for some number\/ $\zeta\in \RR$, we have
\begin{equation}
\label{e:v_tau_n->v*}
  \left\| v(z, t + \tau_n) - U(z + \zeta)
  \right\|_{ C^{2,1}(J\times [0,T]) } \;\longrightarrow\; 0
    \quad\mbox{ as }\, n\to \infty \,.
\end{equation}
Furthermore, we have
\begin{equation*}
  \sup_{ (z,t)\in \RR\times [0,T] }\
  \left| v(z, t + \tau_n) - U(z + \zeta) \right|
    \;\longrightarrow\; 0 \quad\mbox{ as }\, n\to \infty \,.
\end{equation*}
\end{proposition}
\par\vskip 10pt

We combine this proposition with the Lyapunov stability of
travelling waves (Proposition~\ref{prop-Lyapunov})
to prove our main result:

\begin{theorem}\label{thm-Lyapunov}
{\rm (Convergence.)}$\,$
Assume that the reaction function $f\colon \RR\to \RR$ satisfies
all four\/ {\rm Hypotheses}, {\bf (H1)} through {\bf (H4)}, and
the initial data $v_0\colon \RR\to \RR$ satisfy\/
{\rm Hypothesis {\bf (H5)}}.
Assume $-\infty < z_0 < z_1 < +\infty$ and let\/
$J\subset \RR$ be the compact interval specified above.
Let\/
$v\colon \RR\times (0,\infty)\to \RR$
denote the unique bounded classical solution of problem~\eqref{e_c:FKPP}.
Then there is a spatial shift\/ $\zeta = \zeta(v_0)\in \RR$
uniquely determined by the initial data $v_0$, such that\/
\begin{equation}
\label{e:v(z,t)->U(z)}
  \sup_{z\in \RR}\
  | v(z,t) - U(z + \zeta) | \;\longrightarrow\; 0
  \quad\mbox{ as }\, t\to \infty \,.
\end{equation}
Even the following stronger convergence in $C^2(J)$ holds
(cf.~\eqref{e:v_tau_n->v*}),
\begin{equation}
\label{e:v(.,t)->U}
  \| v(\,\cdot\,,t) - U(\,\cdot\, + \zeta)\|_{ C^2(J) }
  \;\longrightarrow\; 0
    \quad\mbox{ as }\, t\to \infty \,.
\end{equation}
\end{theorem}

\par\vskip 10pt
\proof
First, let $\zeta\in \RR$ be the number obtained in
Proposition~\ref{prop-time_conv} above.
Next, we recall from Proposition~\ref{prop-Lyapunov}
that the corresponding TW\--solution
$z\mapsto U(z + \zeta)\colon \RR\to [0,1]$ is Lyapunov\--stable, i.e.,
\eqref{e:v_0-U} $\Longrightarrow$ \eqref{e:v(t)-U}.
We now apply the conclusion of Proposition~\ref{prop-time_conv}
as follows.
For every $\eps\in (0,\eps_0]$ there exists an index
$N = N(\eps)\in \NN = \{ 1,2,3,\dots\}$ such that,
by \eqref{e:v_tau_n->v*},
the following form of \eqref{e:v_0-U} is valid:
\begin{equation}
\label{e:v(tau_n)-U}
\begin{aligned}
  0\leq v(z,\tau_n)\leq 1 \quad\mbox{ and }\quad
  |v(z,\tau_n) - U(z + \zeta)| \leq \eps
\\
  \quad\mbox{ for all }\, z\in \RR \,\mbox{ and for every }\, n\geq N \,.
\end{aligned}
\end{equation}
Thus, by \eqref{e:v(t)-U}, we have
\begin{equation}
\label{e:v(t+tau_n)-U}
  \sup_{z\in \RR}\
  |v(z, t + \tau_n) - U(z + \zeta)| \leq \varrho(\eps)
  \quad\mbox{ for all }\, t\in \RR_+ \,\mbox{ and }\, n\geq N \,.
\end{equation}
Consequently, given any $\eps\in (0,\eps_0]$, there exists some
$t_0 = t_0(\eps)\in [\tau_0,\infty)$, say, $t_0 = \tau_N$, such that
\begin{equation*}
  \sup_{z\in \RR}\
  |v(z,t) - U(z + \zeta)| \leq \varrho(\eps)
  \quad\mbox{ is valid for every }\, t\geq t_0 \,.
\end{equation*}
Recalling $\varrho(\eps)\to 0$ as $\eps\to 0+$,
we infer from this convergence result that the spatial shift
$\zeta\in \RR$ is independent from the choice of the sequence
$\{\tau_n\}_{n=1}^{\infty} \subset [\tau_0,\infty)$
in Proposition~\ref{prop-time_conv}.
It depends solely on the initial data $v_0$.
Moreover, we must have even \eqref{e:v(z,t)->U(z)},
thanks to inequalities
\eqref{e:v(tau_n)-U} and \eqref{e:v(t+tau_n)-U}.
The desired result \eqref{e:v(z,t)->U(z)} follows.

Finally, we apply the same regularity arguments as in the proofs of
Lemma~\ref{lem-time_conv} and Proposition~\ref{prop-time_conv}
to derive \eqref{e:v(.,t)->U} from \eqref{e:v(z,t)->U(z)}.

The theorem is proved.
\qed
\par\vskip 10pt

\subsection*{Acknowledgments}
\begin{small}
The work of Pavel Dr\'abek was supported in part by
the Grant Agency of the Czech Republic (GA\v{C}R)
under Grant {\#}$13-00863$S, and
the work of Peter Tak\'a\v{c} by
the Deutsche Forschungs\-gemeinschaft (DFG, Germany)
under Grants {\#} TA~213/15--1 and {\#} TA~213/16--1.
Both authors were partially supported also by
a joint exchange program between the Czech Republic and Germany;
by the Ministry of Education, Youth, and Sports of the Czech Republic
under the grant No.\ 7AMB14DE005 (exchange program ``MOBILITY'')
and by the Federal Ministry of Education and Research of Germany
under grant No.\ 57063847 (D.A.A.D.\ Program ``PPP'').
Both authors would like to express their sincere thanks
to Professor Hiroshi Matano (University of Tokyo, Japan)
for suggesting to them the monotonicity of travelling waves
in Proposition~\ref{prop-monot_TW}.
\end{small}
%


%
%
\makeatletter \renewcommand{\@biblabel}[1]{\hfill#1.} \makeatother
%
%


\begin{thebibliography}{99}

\bibitem{AronWein}
D.~G. Aronson and H.~F. Weinberger,
{\it Multi\-dimensional nonlinear diffusion arising in
     population genetics},
     Advances in Math., {\bf 30} (1978), 33--76.

\bibitem{Ball}
J.~M. Ball,
{\it Strongly continuous semigroups, weak solutions, and
     the variation\--of\--constants formula},
     Proc. Amer. Math. Soc., {\bf 63}(2) (1977), 70--73.

\color{green}





\color{black}

\bibitem{Deimling}
Klaus Deimling,
{\sl ``Nonlinear Functional Analysis''},
     Springer\--Verlag, Berlin\--Heidelberg, 1985.



\bibitem{DrabManTak}
P. Dr\'abek, R.~F. Man\'asevich, and P. Tak\'a\v{c},
{\it Manifolds of critical points in a quasi\-linear model
     for phase transitions}.
     In
     D.~Bonheure, M.~Cuesta, E.~J.\ Lami Dozo, P.~Tak\'a\v{c},
     J.~Van Schaftingen, and M.~Willem; eds.,
{\sl ``Nonlinear Elliptic Partial Differential Equations''},
     Proceedings of the 2009
     ``International Workshop in Nonlinear Elliptic PDEs,''
     A celebration of Jean\--Pierre Gossez's $65$-th birthday,
     September 2--4, 2009, Brussels, Belgium.
     Contemporary Mathematics, Vol. {\bf 540}, pp.~95--134,
     American Mathematical Society, Providence, R.I., U.S.A., 2011.

\bibitem{DrabMil}
P. Dr\'abek and J. Milota,
{\sl ``Methods of Nonlinear Analysis''}, $2^\mathrm{nd}$ ed.,
     in {\em Birkh\"auser Advanced Texts}.
     Springer\--Verlag, Basel\--Heidelberg\--New York, 2013.

\bibitem{DrabTak-1}
P. Dr\'abek and P. Tak\'a\v{c},
{\it New patterns of travelling waves in
     the generalized Fisher\--Kolmogorov equation},
     Nonlinear Differ. Equ. Appl. (NoDEA), {\bf 23}(2) (2016),
     Article~$7$ (online).
\color{blue}
     {\it Online}: http://dx.doi.org/10.1007/s00030-016-0365-2.
\color{black}



\bibitem{Fife-McLeod}
P.~C. Fife and J.~B. Mc{L}eod,
{\it The approach of solutions of nonlinear diffusion equations to
     travelling front solutions},
     Arch. Rational Mech. Anal., {\bf 65}(4) (1977), 335--361.

\bibitem{Fisher}
R.~A. Fisher,
{\it The advance of advantageous genes},
     Ann. of Eugenics, {\bf 7} (1937), 355--369.


\bibitem{Hale}
Hale, J.~K.,
{\sl ``Asymptotic Behavior of Dissipative Systems''},
Math. Surveys and Monographs, Vol.~\textbf{25}.
American Math. Soc., Providence, R.I., 1988.

\bibitem{Hamel-Nadira}
F. Hamel and N. Nadirashvili,
{\it Travelling fronts and entire solutions of
     the Fisher\--KPP equation in $\mathbb{R}^N$},
     Arch. Rational Mech. Anal., {\bf 157} (2001), 91--163.


\bibitem{Henry}
D.~Henry,
{\sl ``Geometric Theory of Semilinear Parabolic Equations''},
     in {\em Lect. Notes in Math.}, Vol.~\textbf{840}.
     Springer\--Verlag, Berlin\--Heidelberg\--New York, 1981.



\bibitem{KPP}
A. Kolmogorov, I. Petrovski, and N. Piscounov,
{\it \`Etude de l'\'equation de la diffusion avec croissance
     de la quantit\'e de la mati\`ere et son application
     \`a un probl\`eme biologique},
     Bull. Univ. Moskou Ser. Internat. Sec.~A, {\bf 1} (1937), 1--25.

\bibitem{Ladyz-parab}
O.~A. Ladyzhenskaya, N.~N. Ural'tseva, and V.~A. Solonnikov,
{\sl ``Linear and Quasi\--linear Equations of Parabolic Type''}.
In {\em Transl. Mathematical Monographs}, Vol.~{\bf 23},
Amer. Math. Soc., Providence, R.I., 1968.

\color{green}
\color{black}

\bibitem{Matano-Og-1}
H. Matano and T. Ogiwara,
{\it Stability in order\--preserving systems
     in the presence of symmetry},
     Proc. Royal Soc. Edinburgh, {\bf 129\,A} (1999), 395--438.

\bibitem{Matano-Og-2}
H. Matano and T. Ogiwara,
{\it Monotonicity and convergence results
     in order\--preserving systems in the presence of symmetry},
     Discrete and Continuous Dynamical Systems,
     {\bf 5}(1) (1999), 1--34.

\bibitem{Murray}
J.~D. Murray,
{\sl ``Mathematical Biology''},
in {\em Biomathematics Texts}, Vol.~{\bf 19},
Springer-Verlag, Berlin--Hei\-del\-berg\--New York, 1993.

\bibitem{Murray-I}
J.~D. Murray,
{\sl ``Mathematical Biology I: An Introduction''}, $3$-rd Ed.
In {\em Interdisciplinary Applied Mathematics}, Vol.~{\bf 17},
Springer-Verlag, Berlin--Hei\-del\-berg\--New York, 2002.

\bibitem{Pazy}
A. Pazy,
{\sl ``Semigroups of Linear Operators and Applications to
     Partial Differential Equations''},
in {\em Applied Mathematical Sciences}, Vol.~{\bf 44},
Springer-Verlag, New York\--Berlin--Hei\-del\-berg, 1983.

\bibitem{Tsoular_Wall}
A. Tsoularis and J. Wallace,
{\it Analysis of logistics growth models},
     Math. Biosciences, {\bf 179}(1) (2002), 21--55.

\color{green}

\color{black}






\end{thebibliography}
\end{document}